\documentclass[11pt]{article}

\usepackage{latexsym}
\usepackage[dvips]{graphicx}
\usepackage{float}
\usepackage{amsmath}
\usepackage{amssymb}
\usepackage{epic}
\usepackage{color}
\usepackage{setspace}
\usepackage{lscape}
\usepackage{arydshln}
\usepackage{fancybox}   
\usepackage{mathtools}
\usepackage{amsthm}

\begin{document}
\bibliographystyle{plain}
\floatplacement{table}{H}
\newtheorem{definition}{Definition}[section]
\newtheorem{lemma}{Lemma}[section]
\newtheorem{theorem}{Theorem}[section]
\newtheorem{corollary}{Corollary}[section]
\newtheorem{proposition}{Proposition}[section]

\newcommand{\sni}{\sum_{i=1}^{n}}
\newcommand{\snj}{\sum_{j=1}^{n}}
\newcommand{\smj}{\sum_{j=1}^{m}}
\newcommand{\sumjm}{\sum_{j=1}^{m}}
\newcommand{\bdis}{\begin{displaymath}}
\newcommand{\edis}{\end{displaymath}}
\newcommand{\beq}{\begin{equation}}
\newcommand{\eeq}{\end{equation}}
\newcommand{\beqn}{\begin{eqnarray}}
\newcommand{\eeqn}{\end{eqnarray}}
\newcommand{\defeq}{\stackrel{\triangle}{=}}
\newcommand{\simleq}{\stackrel{<}{\sim}}
\newcommand{\sep}{\;\;\;\;\;\; ; \;\;\;\;\;\;}
\newcommand{\real}{\mbox{$ I \hskip -4.0pt R $}}
\newcommand{\complex}{\mbox{$ I \hskip -6.8pt C $}}
\newcommand{\integ}{\mbox{$ Z $}}
\newcommand{\realn}{\real ^{n}}
\newcommand{\sqrn}{\sqrt{n}}
\newcommand{\sqrtwo}{\sqrt{2}}
\newcommand{\prf}{{\bf Proof. }}

\newcommand{\onehlf}{\frac{1}{2}}
\newcommand{\thrhlf}{\frac{3}{2}}
\newcommand{\fivhlf}{\frac{5}{2}}
\newcommand{\onethd}{\frac{1}{3}}
\newcommand{\lb}{\left ( }
\newcommand{\lcb}{\left \{ }
\newcommand{\lsb}{\left [ }
\newcommand{\labs}{\left | }
\newcommand{\rb}{\right ) }
\newcommand{\rcb}{\right \} }
\newcommand{\rsb}{\right ] }
\newcommand{\rabs}{\right | }
\newcommand{\lnm}{\left \| }
\newcommand{\rnm}{\right \| }
\newcommand{\lambdab}{\bar{\lambda}}
%
%
\newcommand{\xj}{x_{j}}
\newcommand{\xzero}{x_{0}}
\newcommand{\xbar}{\bar{x}}
\newcommand{\xstar}{x^{*}}
\newcommand{\tzero}{t_{0}}
\newcommand{\tstar}{t^{*}}

\newcommand{\Azero}{A_{0}}
\newcommand{\Aone}{A_{1}}
\newcommand{\Atwo}{A_{2}}
\newcommand{\Ath}{A_{3}}
\newcommand{\Afr}{A_{4}}
\newcommand{\Afv}{A_{5}}
\newcommand{\Asx}{A_{6}}
\newcommand{\Anmo}{A_{n-1}}
\newcommand{\Anmt}{A_{n-2}}
\newcommand{\An}{A_{n}}
\newcommand{\Aj}{A_{j}}

\newcommand{\azero}{a_{0}}
\newcommand{\aone}{a_{1}}
\newcommand{\atwo}{a_{2}}
\newcommand{\ath}{a_{3}}
\newcommand{\afr}{a_{4}}
\newcommand{\afv}{a_{5}}
\newcommand{\asx}{a_{6}}
\newcommand{\anmo}{a_{n-1}}
\newcommand{\anmt}{a_{n-2}}
\newcommand{\anmth}{a_{n-3}}
\newcommand{\an}{a_{n}}
\newcommand{\aj}{a_{j}}
\newcommand{\ajpt}{a_{j+2}}
\newcommand{\ajpthr}{a_{j+3}}
\newcommand{\ajmo}{a_{j-1}}
\newcommand{\ai}{a_{i}}
\newcommand{\ajpo}{a_{j+1}}
\newcommand{\ajmt}{a_{j-2}}
\newcommand{\ajmthr}{a_{j-3}}
\newcommand{\ak}{a_{k}}
\newcommand{\akpo}{a_{k+1}}
\newcommand{\ajmk}{a_{j-k}}

\newcommand{\bzero}{b_{0}}
\newcommand{\bone}{b_{1}}
\newcommand{\btwo}{b_{2}}
\newcommand{\bth}{b_{3}}
\newcommand{\bfr}{b_{4}}
\newcommand{\bfv}{b_{5}}
\newcommand{\bsx}{b_{6}}
\newcommand{\bnmo}{b_{n-1}}
\newcommand{\bnmt}{b_{n-2}}
\newcommand{\bn}{b_{n}}
\newcommand{\bj}{b_{j}}
\newcommand{\bi}{b_{i}}
\newcommand{\bjpo}{b_{j+1}}
\newcommand{\bk}{b_{k}}
\newcommand{\bkpo}{b_{k+1}}

\newcommand{\alfazero}{\alpha_{0}}
\newcommand{\alfaone}{\alpha_{1}}
\newcommand{\alfatwo}{\alpha_{2}}
\newcommand{\alfathr}{\alpha_{3}}
\newcommand{\alfath}{\alpha_{3}}
\newcommand{\alfafr}{\alpha_{4}}
\newcommand{\alfafv}{\alpha_{5}}
\newcommand{\alfasx}{\alpha_{6}}
\newcommand{\alfanmo}{\alpha_{n-1}}
\newcommand{\alfanmt}{\alpha_{n-2}}
\newcommand{\alfan}{\alpha_{n}}
\newcommand{\alfaj}{\alpha_{j}}
\newcommand{\alfai}{\alpha_{i}}
\newcommand{\alfak}{\alpha_{k}}
\newcommand{\alfajpo}{\alpha_{j+1}}
\newcommand{\alfakpo}{\alpha_{k+1}}

\newcommand{\gamzero}{\gamma_{0}}
\newcommand{\gamone}{\gamma_{1}}
\newcommand{\gamtwo}{\gamma_{2}}
\newcommand{\gamth}{\gamma_{3}}
\newcommand{\gamthr}{\gamma_{3}}
\newcommand{\gamfr}{\gamma_{4}}
\newcommand{\gamfv}{\gamma_{5}}
\newcommand{\gamsx}{\gamma_{6}}
\newcommand{\gamnmo}{\gamma_{n-1}}
\newcommand{\gamnmt}{\gamma_{n-2}}
\newcommand{\gamn}{\gamma_{n}}
\newcommand{\gamj}{\gamma_{j}}
\newcommand{\gamk}{\gamma_{k}}
\newcommand{\gamjpo}{\gamma_{j+1}}
\newcommand{\gamkpo}{\gamma_{k+1}}

\newcommand{\delzero}{\delta_{0}}
\newcommand{\delone}{\delta_{1}}
\newcommand{\deltwo}{\delta_{2}}
\newcommand{\delth}{\delta_{3}}
\newcommand{\delthr}{\delta_{3}}
\newcommand{\delfr}{\delta_{4}}
\newcommand{\delfv}{\delta_{5}}
\newcommand{\delsx}{\delta_{6}}
\newcommand{\delnmo}{\delta_{n-1}}
\newcommand{\delnmt}{\delta_{n-2}}
\newcommand{\deln}{\delta_{n}}
\newcommand{\delj}{\delta_{j}}
\newcommand{\delk}{\delta_{k}}
\newcommand{\deljpo}{\delta_{j+1}}
\newcommand{\delkpo}{\delta_{k+1}}

\newcommand{\muone}{\mu_{1}}
\newcommand{\mutwo}{\mu_{2}}

\newcommand{\Bzero}{B_{0}}
\newcommand{\Bone}{B_{1}}
\newcommand{\Btwo}{B_{2}}
\newcommand{\Bth}{B_{3}}
\newcommand{\Bfr}{B_{4}}
\newcommand{\Bfv}{B_{5}}
\newcommand{\Bsx}{B_{6}}
\newcommand{\Bnmo}{B_{n-1}}
\newcommand{\Bndtwomo}{B_{n/2-1}}
\newcommand{\Bnmt}{B_{n-2}}
\newcommand{\Bn}{B_{n}}
\newcommand{\Bj}{B_{j}}

\newcommand{\Anmk}{A_{n-k}}
\newcommand{\Bnmk}{A_{n-k}}
\newcommand{\anmk}{a_{n-k}}
\newcommand{\zj}{z^{j}}          
\newcommand{\zk}{z^{k}}          

\newcommand{\matrixspace}{\;\;}
\newcommand{\ellone}{\ell_{1}}
\newcommand{\elltwo}{\ell_{2}}
\newcommand{\cmm}{\complex^{m \times m}}
\newcommand{\opdsk}{\mathcal{O}}

\newcommand{\xone}{x_{1}}
\newcommand{\xtwo}{x_{2}}
\newcommand{\xthr}{x_{3}}
\newcommand{\xfor}{x_{4}}
\newcommand{\xfiv}{x_{5}}

\newcommand{\lnorm}{\left \|}
\newcommand{\rnorm}{\right \|}
\newcommand{\lnrm}{\biggl | \biggl |}
\newcommand{\rnrm}{\biggr |\biggr |}

\newcommand{\recipp}{p^{\protect \#}}
\newcommand{\recipP}{P^{\protect \#}}

\newcommand{\xkova}{\hat{x}}
\newcommand{\enestrom}{Enestr\"{o}m }
\newcommand{\eneskak}{Enestr\"{o}m-Kakeya}
\newcommand{\eneskakk}{Enestr\"{o}m-Kakeya }

\begin{center}
\large
{\bf A UNIFYING FRAMEWORK FOR GENERALIZATIONS OF THE ENESTR\"{O}M-KAKEYA THEOREM}
\vskip 0.5cm
\normalsize
A. Melman \\
Department of Applied Mathematics \\
School of Engineering, Santa Clara University  \\
Santa Clara, CA 95053  \\
e-mail : amelman@scu.edu \\
\vskip 0.5cm
\end{center}

\begin{abstract}
The classical \eneskakk theorem establishes upper and lower bounds on the zeros of a polynomial with positive coefficients that are explicit 
functions of those coefficients. We establish a unifying framework that incorporates this theorem and several similar ones as special cases,
while generating new theorems of a similar type. These establish zero inclusion and exclusion regions consisting of a single disk or the union 
of several disks in the complex plane. Our framework is built on two basic tools, namely a generalization of an observation by Cauchy, and a 
family of polynomial multipliers. Its approach is transparent and reduces algebraic manipulations to a minimum.
\vskip 0.15cm
{\bf AMS(MOS) subject classification :} 12D10, 30C15, 65H05
\end{abstract}

%
%
%
%
\section{introduction}

The \eneskakk theorem (\cite{Enestrom1}, \cite{Enestrom2}, and \cite{Kakeya}) establishes upper and lower bounds on the moduli of the zeros of a polynomial
with positive coefficients that are simple explicit functions of the coefficients.
It has been extended and generalized in different ways, e.g., by relaxing the restrictions on the 
coefficients, or by deriving different regions of the complex plane that also contain all the zeros of the given polynomial. 
Sharpness of the bounds was considered in~\cite{ASV1}, \cite{ASV2}, and~\cite{Hurwitz}.
A good overview of these generalizations can be found in~\cite{GG} and for additional historical remarks about this pretty theorem 
we also refer to~\cite[p. 271-272]{RS}.

Here we establish a framework that generates theorems similar to that of Enestr\"{o}m and Kakeya, namely, theorems that derive regions in the complex plane
that contain or do not contain the zeros of polynomials with positive coefficients. These regions, which consist of a single disk or of the union of several 
disks, are explicitly determined, i.e., they do not require numerical methods. 
We obtain our results by using the same two basic tools for all of them: a family of polynomial multipliers and a generalization of
an observation by Cauchy. This transparent approach reduces algebraic manipulations to a minimum, unifies the derivation of these results, and generates
new ones while incorporating the \eneskakk theorem and some of its variants as special cases. Some of the regions derived in this way consist
of several disks, which can provide additional information about the zeros. Regions of this kind have not been considered in any of the existing generalizations 
of the \eneskakk theorem.

The paper is organized as follows. We begin in Section~\ref{preliminaries} by collecting a few results that are needed further on. In Section~\ref{onedisk_origin},
we derive zero inclusion and exclusion regions composed of a single disk centered at the origin, while in Section~\ref{onedisk_notorigin} we obtain disks that
are not centered at the origin. In Sections~\ref{twodisks} and \ref{threedisks}, we derive zero inclusion regions consisting of two and three disks, respectively.
The appendix contains a technical lemma needed in Section~\ref{onedisk_notorigin}.

%
%
%
%
\section{Preliminaries}   
\label{preliminaries}

In this section, we collect a few theorems and definitions that will be needed further on.
Throughout this section, we consider a polynomial $p(z) = \an z^{n}+\anmo z^{n-1}+\dots+\aone z + \azero$
with complex coefficients.
\begin{definition}
\label{def_sk}
The \emph{Cauchy radius of the $k$th kind} of $p$ is defined as the unique positive solution $s_{k}$ of 
\beq
\label{sk_eq}
|\an| z^{n}+ |\anmo| z^{n-1}+\dots+|a_{n-k+1}|z^{n-k+1} -|\anmk| z^{n-k} - \dots -|\aone |z - |\azero| = 0 \; .
\eeq   
When $k=1$, $s_{1}$ is simply called the \emph{Cauchy radius} of $p$.
\end{definition}
The following theorem is Theorem~3.2 in~\cite{Melman_RMJM}, where $k$ here corresponds to $n-k$ in that theorem.

%
%
%
%
\begin{theorem}
\label{theorem_bakchik_cauchy}
All the zeros of the complex polynomial $p(z) = \sum_{j=0}^{n} a_{j} z^{j}$ lie in the sets 
\begin{multline}
\nonumber
\Gamma_{1}(k) = \lcb z \in \complex : |\an z^{k} + \anmo z^{k-1} + \dots + a_{n-k+1}z| \right . \\
                                                                        \left .          \leq |\an| s_{k}^{k} + |\anmo| s_{k}^{k-1} + \dots + |a_{n-k+1}|s_{k} \rcb 
\end{multline}
and
\begin{multline}
\nonumber
\Gamma_{2}(k) = \lcb z \in \complex : |\an z^{k} + \anmo z^{k-1} + \dots + a_{n-k+1}z + \anmk | \right .      \\
                                         \left .        \leq |\an| s_{k+1}^{k} + |\anmo| s_{k+1}^{k-1} + \dots + |a_{n-k+1}|s_{k+1} + |\anmk|  \rcb \; , 
\end{multline}
whose boundaries are lemniscates and where $s_{j}$ is the Cauchy radius of the $j$th kind, defined in Definition~\ref{def_sk}.

If $\Gamma_{1}(k)$ or $\Gamma_{2}(k)$ consists of disjoint regions whose boundaries are 
simple closed (Jordan) curves and $\ell$ is the number of foci of its bounding lemniscate contained in any 
such region, then that region contains $\ell$ zeros of $p$ when it does not contain the origin, and  
$\ell + n - k$ zeros of $p$ when it does contain the origin.
\end{theorem}

For the special case $\Gamma_{1}(1)$ in Theorem~\ref{theorem_bakchik_cauchy}, we obtain that all the zeros of $p$ lie in the disk defined by 
$|z| \leq s_{1}$, where $s_{1}$ is the unique positive solution of 
\beq 
\label{Cauchy_eq}
|\an | z^{n} - |a_{n-1}| z^{n-1} - \dots  - |\aone| z - |\azero| = 0  \; . 
\eeq
This is a classical observation by Cauchy from~1829 (\cite{Cauchy}, see also, e.g., 
\cite[Th.(27,1), p.122 and Exercise 1, p.126]{Marden}, \cite[Theorem 3.1.1]{MMR}, \cite[Theorem 8.1.3]{RS}).

In the special case of $\Gamma_{2}(1)$ in Theorem~\ref{theorem_bakchik_cauchy}, we obtain that all the zeros of $p$ lie in the disk defined 
by $|z+\anmo/\an| \leq s_{2} + |\anmo/\an |$, where $s_{2}$ is the unique positive solution of 
\beq 
\label{Twin_eq}
|\an | z^{n} + |a_{n-1}| z^{n-1} - |\anmt | z^{n-2} - \dots  - |\aone| z - |\azero| = 0  \; . 
\eeq
This is Theorem~1 in~\cite{Melman_Twin}. 

The boundaries of $\Gamma_{1}(k)$ and $\Gamma_{2}(k)$ are lemniscates that, as $k$ increases, become too complicated to use for deriving explicit zero
inclusion regions. Instead, we will
approximate a region $\Gamma$, bounded by a lemniscate of the form $|q(z)|=R$, where $q(z)=z^{m} + b_{m-1}z^{m-1} + \dots + b_{0}$, as follows. 
Denoting the zeros of $q$ by $c_{j}$, we have that  
\bdis
|q(z)| = |z^{m} + b_{m-1}z^{m-1} + \dots + b_{0}| = |z-c_{1}||z-c_{2}|\dots |z-c_{m}| \; ,
\edis
so that 
\bdis
\Gamma = \lcb z \in \complex : |q(z)| \leq R \rcb = \lcb z \in \complex :  |z-c_{1}||z-c_{2}|\dots |z-c_{m}| \leq R \rcb \; .
\edis
This means that 
\bdis
\Gamma \subseteq \bigcup_{j=1}^{m} \lcb z \in \complex : |z-c_{j}| \leq R^{1/m}  \rcb \; .
\edis
Although larger than $\Gamma$, this union of disks is easier to work with, and, as we shall see, can still be useful.
It may sometimes be better to allow different disks to have different radii, but in the interest of simplicity we will
use the same radius for all disks. 
 
Additional zero inclusion regions can be obtained by applying Theorem~\ref{theorem_bakchik_cauchy} to
the reverse polynomial $\recipp(\zeta)= \azero \zeta^{n} + \aone \zeta^{n-1} + \dots + \anmo \zeta + \an$ (with $\azero \neq 0$),
whose zeros are the reciprocals of the zeros of $p$. 

%
%
%
%
\section{Single disk centered at the origin}
\label{onedisk_origin}           

In this section we derive inclusion and exclusion regions consisting of a single disk centered at the origin for polynomials with positive coefficients. 
We begin by stating the \eneskakk theorem (\cite{Enestrom1}, \cite{Enestrom2}, and \cite{Kakeya}) and, for completeness, also provide a standard proof.
%
%
\begin{theorem}(\eneskakk )
\label{Theorem_EK}
All the zeros of the real polynomial $p(z) = \sum_{j=0}^{n} a_{j} z^{j}$ with positive coefficients lie in the annulus defined by
\beq
\nonumber
\lcb z \in \complex : \min_{0 \leq j \leq n-1} \dfrac{\aj}{\ajpo} \leq |z| \leq \max_{0 \leq j \leq n-1} \dfrac{\aj}{\ajpo} \rcb .
\eeq
\end{theorem}
\prf Consider $(z-\gamma)p(z)$, where $\gamma \in \real$:
\beq
\label{multipliereq1}
(z-\gamma)p(z) = \an z^{n+1} + \sum_{j=1}^{n} (\ajmo - \gamma \aj ) z^{j} -\gamma \azero \; . 
\eeq
Clearly, any upper bound on the zeros of $(z-\gamma)p(z)$ will also be an upper bound on the zeros of $p$. From~(\ref{multipliereq1}) we observe that
all the coefficients of $(z-\gamma)p(z)$, except that of $z^{n+1}$, will be nonpositive if 
\bdis
\gamma = \max_{0 \leq j \leq n-1} \dfrac{\aj}{\ajpo} \; ,
\edis
while the constant term is negative.
For that value of $\gamma$, the Cauchy radius of $(z-\gamma)p(z)$ is the unique positive solution of $(z-\gamma)p(z)~=~0$ and, since $\gamma > 0$ is a positive
zero of $(z-\gamma)p(z)$, it must be equal to the Cauchy radius. This proves the upper bound. The lower bound is proved analogously by considering the reverse 
polynomial of $p$. \qed

The key ingredients in this proof are the multiplier $z-\gamma$ and the Cauchy radius, a pattern that will repeat itself in subsequent results, 
albeit with different multipliers and Cauchy radii of a higher kind. 
The use of multipliers is standard practice throughout the literature related to the \eneskakk theorem.

It is worth pointing out that the \eneskakk bounds are not necessarily better than those obtained from the Cauchy radii of $p$ and $\recipp$, 
although their obvious advantage 
is that they are explicit and therefore do not require the solution of a polynomial equation. Consider the following example.
\vskip 0.25cm \noindent {\bf Example.}
Define the polynomials  
$p_{1}(z)=2z^{5}+z^{4}+4z^{3}+z^{2}+2z+3$ and $p_{2}(z)=z^{5}+z^{4}+2z^{3}+3z^{2}+2z+1$. 
For $p_{1}$, whose smallest and largest zeros have
magnitudes $0.7740$ and $1.3921$, respectively, the Cauchy radii determine an interval of $[0.6288,1.9242]$ for the magnitudes of the zeros, while the \eneskakk
theorem gives $[0.25,4.00]$, which is worse. 

On the other hand, for $p_{2}$, with smallest and largest zero magnitudes $0.7136$ and $1.4013$, 
respectively, we obtain $[0.3143,2.4654]$ and $[0.5,2.0]$ from the Cauchy radii and Enestr\"{o}m-Kakeya, respectively. Here, the \eneskakk bounds are better.
%
\vskip 0.25cm
There is one clear case where the \eneskakk theorem is guaranteed to produce a better upper bound than the Cauchy radius, namely, when
$\anmo/\an=\max_{1 \leq j \leq n} a_{j-1}/a_{j}$. 
This is an immediate consequence of Theorem~8.3.1 in~\cite{RS}, which states that, if the Cauchy radius of $p$ is not a zero of $p$ (which is obviously the case
since all of its coefficients are positive), 
then the Cauchy radius of $(z-\anmo/\an)p(z)$, which here is equal to $\anmo/\an$, is strictly smaller than that of $p$ for any complex polynomial $p$, 
whose coefficients are all nonzero. This can also easily be seen by substituting $\anmo/\an$ in the left-hand side of~(\ref{Cauchy_eq}), which yields a 
negative value, indicating that $\anmo/\an$ is less than the Cauchy radius.
Analogously, the \eneskakk theorem will yield a better lower bound than that obtained from the Cauchy radius 
if $\azero/\aone=\min_{1 \leq j \leq n} a_{j-1}/a_{j}$.

We can obtain more results like the \eneskakk theorem by simply changing the multiplier used in its proof, while using Theorem~\ref{theorem_bakchik_cauchy}
to obtain bounds. The following two theorems do just that. 

%
%
\begin{theorem}
\label{Theorem_Bakchik1}
Let the real polynomial $p(z) = \sum_{j=0}^{n} a_{j} z^{j}$ with $n \geq 2$ have positive coefficients. 
\newline {\bf (a)}
Let 
\beq
\nonumber
\gamone = \dfrac{\anmo}{\an} 
\;\; \text{and} \;\;
\gamzero = \max \lcb 0 , \max_{0 \leq j \leq n-2} \dfrac{\aj-\gamone\ajpo}{\ajpt} \rcb \; ,
\eeq
then all the zeros of $p$ lie in the disk defined by $|z| \leq \frac{1}{2} \lb \gamone + \lb \gamone^{2} + 4\gamzero \rb^{1/2} \rb$.
\newline {\bf (b)}
Let 
\beq
\nonumber
\delone = \dfrac{\aone}{\azero} 
\;\; \text{and} \;\;
\delzero = \max \lcb 0 , \max_{0 \leq j \leq n-2} \dfrac{\ajpt-\delone\ajpo}{\aj} \rcb \; ,
\eeq
then all the zeros of $p$ are excluded from an open disk defined by 
\bdis
|z| < \dfrac{2}{\delone + \lb \delone^{2} + 4\delzero \rb^{1/2}} \; .
\edis
\end{theorem}
\prf Consider $q(z)=(z^{2}-\gamone z-\gamzero)p(z)$ for $\gamzero,\gamone \in \real$:
\beq
\label{multipliereq2}
q(z) = \an z^{n+2} + (\anmo -\gamone \an) z^{n+1} + \sum_{j=2}^{n} (\ajmt - \gamone \ajmo - \gamzero \aj) z^{j} 
-\gamone \azero z -\gamzero \azero \; . 
\eeq
Any set containing all the zeros of $q$ will also contain all those of $p$. If all coefficients in equation~(\ref{multipliereq2}),
other than the leading one, are nonpositive, then the Cauchy radius of $q$ is its own unique positive zero.
With $\gamone=\anmo/\an$, the smallest such value for $\gamzero$ is given by 
$\gamzero = \max_{0 \leq j \leq n-2} \lb \ajmt-\gamone\ajmo \rb/a_{j}$. If this value is negative, then 
we must set $\gamzero=0$ to make the constant coefficient in~(\ref{multipliereq2}) nonpositive. 
The Cauchy radius is then the unique positive zero of the quadratic multiplier, which is
$\frac{1}{2} \lb \gamone + \lb \gamone^{2} + 4\gamzero \rb^{1/2} \rb$. This proves part (a).

For part (b), we apply the result in part (a) to the reverse polynomial $\recipp(\zeta)= \azero \zeta^{n} + \aone \zeta^{n-1} + \dots + \anmo \zeta + \an$,
whose zeros are the reciprocals of those of $p$. 
It is a straightforward exercise to
replace $\aj$ by $a_{n-j}$ in part (a) to obtain the equivalent expressions for the reciprocals $\zeta=1/z$ of the zeros $z$ of $p$, bearing in mind
that $|\zeta | \leq M$ implies that $|z| \geq 1/M$. This concludes the proof. \qed


Theorem~\ref{Theorem_Bakchik1} delivers a better bound than the Cauchy radius when $\gamzero=(\anmt-\gamone\anmo)/\an$. This is a consequence of Theorem~3.1
in~\cite{Melman_PAMS3}, which states that, if the Cauchy radius of $p$ is not a zero of $p$,  
then the Cauchy radius of $(z^{2}-(\anmo/\an)z-\anmt/\an+\anmo^{2}/\an^{2})p(z)$ is strictly smaller than that of $p$.
Analogously, we find that Theorem~\ref{Theorem_Bakchik1} yields a better lower bound than can be obtained from the Cauchy radius when 
$\delzero=(\atwo-\delone\aone)/\azero$. We also note that, if $\anmo/\an = \max_{0 \leq j \leq n-1} \aj/\ajpo$, then $\gamzero=0$ and the upper bound from 
this theorem is identical to the \eneskakk upper bound, which in this case is smaller (better) than the Cauchy radius, as we saw before. An analogous 
argument holds for the lower bound.

Although somewhat dissimilar at first sight, part (a) of this theorem is essentially Theorem~1 in~\cite{AM} with a few differences. 
Our result is formulated in terms of parameters that are the coefficients of a multipier, 
whereas in~\cite{AM} the parameters are the zeros of the multiplier. This is merely a question of preference, although using the coefficients leads to simpler
expressions. In cite~\cite{AM}, the parameters of the multiplier are left to be determined as long as they make the appropriate coefficients nonpositive,
exactly like our parameters here, but with the distinction that 
we actually assign values to the parameters to satisfy that requirement. More importantly, the proof in~\cite{AM} relies on more
complicated arguments than the Cauchy radius. That, combined with using zeros rather than coefficients as parameters, would appear to make it more difficult 
to use higher order multipliers to obtain results as in the next theorem. Finally, although a minor matter, no lower bound was mentioned in~\cite{AM}.

Theorem~\ref{Theorem_Bakchik1} sets the stage for similar results. The following theorem is a prototype
for such results and it illustrates the straightforward manner in which they can be obtained. They are only marginally more complicated than 
Theorem~\ref{Theorem_Bakchik1}, but they tend to produce better results, as we will see later.
%
%
\begin{theorem}
\label{Theorem_Bakchik2}
Let the real polynomial $p(z) = \sum_{j=0}^{n} a_{j} z^{j}$ with $n \geq 3$ have positive coefficients, then the following holds. 
\newline {\bf (a)}
If
\bdis
\gamtwo=\max \lcb \dfrac{\anmo}{\an} , \dfrac{\anmt}{\anmo} \rcb \; , \; 
\gamone= 0 \; , \; 
\gamzero = \max \lcb 0 , \max_{0 \leq j \leq n-3} \dfrac{\aj-\gamtwo\ajpo}{\ajpthr} \rcb \; , 
\edis
or if
\begin{eqnarray*}
& & \gamtwo=\dfrac{\anmo}{\an} \; , \\
& & \gamone=\dfrac{\anmt-\gamtwo\anmo}{\an} \; , \\ 
& & \gamzero = \max \lcb 0 , \dfrac{\gamone\azero}{-\aone} , \dfrac{\gamtwo\azero+\gamone\aone}{-\atwo} , 
\max_{0 \leq j \leq n-3} \dfrac{\aj-\gamtwo\ajpo-\gamone\ajpt}{\ajpthr} \rcb \; , 
\end{eqnarray*}
then all the zeros of $p$ lie in the disk defined by $|z| \leq r_{1}$, where $r_{1}$ is the unique positive zero of $z^{3}-\gamtwo z^{2}-\gamone z -\gamzero$. 
\newline {\bf (b)}
If
\bdis
\deltwo=\max \lcb \dfrac{\aone}{\azero} , \dfrac{\atwo}{\aone} \rcb \; , \; 
\delone= 0 \; , \; 
\delzero = \max \lcb 0 , \max_{0 \leq j \leq n-3} \dfrac{\ajpthr-\deltwo\ajpt}{\aj} \rcb \; , 
\edis
or if
\begin{eqnarray*}
& & \deltwo=\dfrac{\aone}{\azero} \; , \\
& & \delone=\dfrac{\atwo-\deltwo\aone}{\azero} \; ,  \\ 
& & \delzero = \max \lcb 0 , \dfrac{\delone\an}{-\anmo} , \dfrac{\deltwo\an+\delone\anmo}{-\anmt} , 
\max_{0 \leq j \leq n-3} \dfrac{\ajpthr-\deltwo\ajpt-\delone\ajpo}{\aj} \rcb \; ,
\end{eqnarray*}
then all the zeros of $p$ are excluded from an open disk defined by $|z| < 1/r_{2}$, where $r_{2}$ is the unique positive zero of 
$z^{3}-\deltwo z^{2}-\delone z -\delzero$. 
\end{theorem}
\prf
Consider $q(z)=(z^{3}-\gamtwo z^{2}-\gamone z -\gamzero)p(z)$,
where $\gamzero,\gamone,\gamtwo \in \real$:
\begin{multline}
\label{multipliereq3}
q(z) = \an z^{n+3} + (\anmo -\gamtwo \an) z^{n+2} + (\anmt - \gamtwo \anmo - \gamone \an )z^{n+1} \\ 
\hskip 2cm + \sum_{j=3}^{n} \lb \ajmthr - \gamtwo \ajmt - \gamone \ajmo - \gamzero \aj \rb  z^{j} 
- \lb \gamtwo \azero +\gamone \aone +\gamzero\atwo \rb z^{2}  \\
\hskip 3cm - \lb \gamone \azero +\gamzero \aone \rb z -\gamzero\azero \; . 
\end{multline}
If we set 
\bdis
\gamtwo=\max \lcb \dfrac{\anmo}{\an} , \dfrac{\anmt}{\anmo} \rcb \; , \; 
\gamone= 0 \; , \; 
\gamzero = \max \lcb 0 , \max_{0 \leq j \leq n-3} \dfrac{\aj-\gamtwo\ajpo}{\ajpthr} \rcb \; , 
\edis
or if we set
$\; \gamtwo=\anmo/\an$, $\gamone = \lb \anmt-\gamtwo\anmo \rb /\an$, and 
\bdis
\gamzero = \max \lcb 0 , \dfrac{\gamone\azero}{-\aone} , \dfrac{\gamtwo\azero+\gamone\aone}{-\atwo} , 
\max_{0 \leq j \leq n-3} \dfrac{\aj-\gamtwo\ajpo-\gamone\ajpt}{\ajpthr} \rcb \; ,
\edis
then the coefficients of the nonleading powers of $z$ in the right-hand side of~(\ref{multipliereq3}) are all nonpositive,
which means, reasoning as before, that the Cauchy radius of $q$, which is also an upper bound on the moduli of the zeros of $p$, 
is the unique positive zero of $z^{3}-\gamtwo z^{2}-\gamone z -\gamzero$. This proves part~(a). 

As before, to prove part (b), we apply the result in part (a) to the reverse polynomial $\recipp$, and the proof follows analogously.
\qed

Note that, for the first choice of multiplier coefficients in this theorem, $\gamone$ and $\delone$ can be negative, and
for that first choice we also observe that, if 
\bdis
\dfrac{\anmo}{\an} = \max_{0 \leq j \leq n-1} \dfrac{\aj}{\ajpo} \;\;\;\; \text{or} \;\;\;\; \dfrac{\anmt}{\anmo} = \max_{0 \leq j \leq n-1} \dfrac{\aj}{\ajpo} \; ,
\edis
then $\gamzero=0$ and the upper bound becomes equal to the one from the \eneskakk theorem. If $\anmo/\an \geq \anmt/\anmo$,
this bound will be smaller than the Cauchy radius.

For the second choice of parameters, if
\bdis 
\max_{0 \leq j \leq n-3} \dfrac{\aj}{\ajpo} \leq \dfrac{\anmo}{\an} \leq \dfrac{\anmt}{\anmo} \; ,
\edis
then $\gamone \geq 0$, which implies that $\gamzero=0$, and we obtain precisely the result in Theorem~\ref{Theorem_Bakchik1}, which in this case
yields a smaller upper bound than the Cauchy radius, as was explained immediately after the proof of that theorem.
Analogous conclusions can be drawn for the lower bound. 
\vskip 0.25 cm
\noindent
{\bf Example.} We illustrate the theorems in this section with the polynomial 
$p_{3}(z)=z^{6}+4z^{5}+2z^{4}+2z^{3}+3z^{2}+6z+7$. The smallest and largest moduli of its zeros 
are the endpoints of the interval $[1.075,3.554]$. The different theorems produce the following intervals:
\bdis
\begin{array}{ll}
\text{Cauchy radii:}  &  [0.670,4.580]   \; , \\
\text{\eneskakk theorem:}  &  [0.500,4.000]   \; , \\
\text{Theorem~\ref{Theorem_Bakchik1}:}  &  [0.633,4.000] \; ,  \\
\text{Theorem~\ref{Theorem_Bakchik2}(1):}  &  [0.766,4.000] \; ,  \\
\text{Theorem~\ref{Theorem_Bakchik2}(2):}  &  [0.807,3.788] \; .  \\
\end{array}
\edis
The $(1)$ and $(2)$ versions of Theorem~\ref{Theorem_Bakchik2} refer to the first and second choices, respectively, for the parameters in that theorem.
Since here $\anmo/\an=\max_{0 \leq j \leq n-1} \aj/\ajpo$, we find, as predicted, that the \eneskakk upper bound is the same as that obtained from 
theorems~\ref{Theorem_Bakchik1}
and~\ref{Theorem_Bakchik2}(1). There is clearly a positive trend in the results when going from the \eneskakk theorem to Theorem~\ref{Theorem_Bakchik2}(2).

It is, in general, difficult to predict which of the theorems in this section will produce the best results. Theorems obtained with higher order
multipliers do not necessarily outperform those obtained with lower order multipliers, although they frequently do. To gain some insight into their
performance, we have carried out a series of numerical comparisons with two classes of randomly generated polynomials of degrees $10$ and $40$,  
which allows us to observe the effect of increasing the degree. We have included the Cauchy radii in these 
comparisons for reference, since they are typically among the best bounds one can expect. The two classes are as follows.
\begin{itemize}
\item \underline{{\bf Class I:}} Polynomials, whose coefficients are uniformly randomly distributed in $(0,10)$.
\item \underline{{\bf Class II:}} Polynomials, whose coefficients are uniformly randomly distributed in $(1,5)$.
\end{itemize}
The second class of polynomials shows the effect of a more limited range, where the coefficients are bounded away from zero. The results for zero inclusion and 
exclusion sets (or upper and lower bounds) are very similar and we limit ourselves to the former to avoid overloading the tables below. 
We generated $1000$ polynomials for each class and collected the results in the following two tables.
\begin{itemize}
\item Table~\ref{Table_11} lists, for Class I polynomials, the median of the ratios of the upper bound to the modulus of the largest zero of the poynomial, 
i.e., the closer this number is to 1, the better it is.
\item Table~\ref{Table_12} is the analog of Table~\ref{Table_11} for Class II polynomials.
\end{itemize}

For Theorem~\ref{Theorem_Bakchik2}, the first and second choices of the multiplier are designated in the tables by Theorem~\ref{Theorem_Bakchik2}(1)
and Theorem~\ref{Theorem_Bakchik2}(2), respectively. 

%
%
\begin{table}[H]
\begin{center}
\small           
\begin{tabular}{c|c|c|c|c|c}
                  & Cauchy               & Enestr\"{o}m-Kakeya  &  Theorem 3.2         & Theorem 3.3(1)      & Theorem 3.3(2)            \\ \hline        
                  &                      &                      &                      &                     &                           \\
n=10              &        1.465         &        4.570         &        1.600         &       1.479         &        1.3553             \\        
                  &                      &                      &                      &                     &                           \\ 
n=40              &        1.417         &       18.828         &        2.515         &       2.015         &        2.072              \\        
\end{tabular}
\caption{Comparison of upper bounds for Class I polynomials.}
\label{Table_11}
\end{center}
\end{table}
\normalsize

%
%
\begin{table}[H]
\begin{center}
\small           
\begin{tabular}{c|c|c|c|c|c}
                  & Cauchy               & Enestr\"{o}m-Kakeya  &  Theorem 3.2         & Theorem 3.3(1)      & Theorem 3.3(2)            \\ \hline         
                  &                      &                      &                      &                     &                           \\
n=10              &        1.626         &        2.091         &        1.385         &       1.365         &        1.245              \\        
                  &                      &                      &                      &                     &                           \\
n=40              &        1.629         &        2.920         &        1.595         &       1.516         &        1.393              \\        
\end{tabular}
\caption{Comparison of upper bounds for Class II polynomials.}
\label{Table_12}
\end{center}
\end{table}
\normalsize

Summarizing these results, we conclude that theorems of the Enestr\"{o}m-Kakeya type perform much better on Class II-like polynomials, where the 
polynomials have coefficients that are more limited in range and that are bounded away from zero. For the latter, this is not surprising, as the 
expressions for 
the quantities that define the inclusion regions contain the coefficients in their denominators. For such polynomials these theorems frequently deliver 
better results than those 
based on the Cauchy radii. As the degree of the polynomials increases, Cauchy radii tend to become marginally better for Class II polynomials,
while they become dramatically better for Class I polynomials.

Furthermore, for both classes of polynomials, Theorems \ref{Theorem_Bakchik2}(1) and \ref{Theorem_Bakchik2}(2) outperform Theorem~\ref{Theorem_Bakchik1}, while
Theorem \ref{Theorem_Bakchik2}(2) is generally better than Theorem~\ref{Theorem_Bakchik2}(1). 
In all cases, the theorems presented here outperform the Enestr\"{o}m-Kakeya theorem.
As a reminder, we note again that the Cauchy radii do not have explicit expressions and need to be computed by solving a polynomial equation.

\vskip 0.25cm
It is now an entirely straightforward matter to continue the pattern of Theorem~\ref{Theorem_Bakchik1} and Theorem~\ref{Theorem_Bakchik2} with ever higher 
order multipliers. Of course, beyond a quartic, a numerical method is required to compute their zeros. In fact, such a method may be preferable even in the 
cubic or quartic case. A simple method such as Newton's method only requires a few iterations for low order polynomials, but 
whether or not this is practical is, in any event, to be decided elsewhere.
In the following section we derive zero inclusion regions that are not centered at the origin.
 
%
%
\section{Single disk not centered at the origin}
\label{onedisk_notorigin}

So far, we have used polynomial multipliers and the Cauchy radius to derive upper and lower bounds on the moduli of the zeros of a polynomial with positive
coefficients. We now derive additional inclusion regions for the zeros of such polynomials, based, this time, on the Cauchy radius of the second 
kind. The following theorem sets the general tone.
%
%
\begin{theorem}
\label{Theorem_Bakchik3}
Let the real polynomial $p(z) = \sum_{j=0}^{n} a_{j} z^{j}$ with $n \geq 2$ have positive coefficients, and let 
\bdis
\muone = \lb \max_{0 \leq j \leq n-2} \dfrac{\aj}{\ajpt} \rb^{1/2} 
\;\; \text{and} \;\;
\mutwo = \lb \min_{0 \leq j \leq n-2} \dfrac{\aj}{\ajpt} \rb^{1/2}  \; .
\edis
Then all the zeros of $p$ are included in the closed disk
\beq
\nonumber
D_{1} = \lcb z \in \complex : \left | z + \dfrac{\anmo}{\an} \right | \leq \dfrac{\anmo}{\an} + \muone \rcb \; , 
\eeq
but are excluded from the open disk
\beq
\nonumber
D_{2} = \lcb z \in \complex : \left | z - \dfrac{\mutwo^{2}}{\azero/\aone+2\mutwo}  \right | 
< \dfrac{\mutwo \lb \azero/\aone + \mutwo \rb}{\azero/\aone + 2\mutwo} \rcb \; . 
\eeq
\end{theorem}
\prf Consider $q(z)=(z^{2}-\gamma)p(z)$, where $\gamma \in \real$:
\beq
\label{multipliereq4}
q(z) = \an z^{n+2} + \anmo z^{n+1} + \sum_{j=2}^{n} (\ajmt - \gamma \aj ) z^{j} -\gamma \aone z -\gamma \azero \; . 
\eeq
Any set containing all the zeros of $q$ will also contain all those of $p$. Equation~(\ref{multipliereq4}) shows that 
all the coefficients of $q$, except those of $z^{n+2}$ and $z^{n+1}$, will be nonpositive if 
\bdis
\gamma = \max_{0 \leq j \leq n-2} \dfrac{\aj}{\ajpt} \; .
\edis
In that case, the Cauchy radius of the second kind of $q$ is the unique positive solution of $(z^{2}-\gamma)p(z)~=~0$,
which must be $\sqrt{\gamma}$. The set $\Gamma_{2}(1)$ in Theorem~\ref{theorem_bakchik_cauchy} then shows 
that all the zeros of $p$ lie in the disk centered at $-\anmo/\an$ with radius $\anmo/\an+\sqrt{\gamma}$, which is the disk $D_{1}$.

To obtain the zero exclusion disk $D_{2}$, we derive an inclusion disk for the reciprocals of the zeros, which are the zeros of the reverse polynomial
$\recipp(\zeta)= \azero \zeta^{n} + \aone \zeta^{n-1} + \dots + \anmo \zeta + \an$. Proceeding as we did for $D_{1}$, we obtain, by replacing
$\aj$ by $a_{n-j}$, that the reciprocals of the zeros of $p$ lie in the disk centered at $-\aone/\azero$ with radius
\bdis
\lb \max_{0 \leq j \leq n-2} \dfrac{\ajpt}{\aj} \rb^{1/2} = \lb \min_{0 \leq j \leq n-2} \dfrac{\aj}{\ajpt} \rb^{-1/2} = \dfrac{1}{\mutwo} \; ,
\edis
i.e., they lie in the disk
\beq
\nonumber
\lcb \zeta \in \complex : \left | \zeta + \dfrac{\aone}{\azero} \right | \leq \dfrac{\aone}{\azero} + \dfrac{1}{\mutwo} \rcb \; . 
\eeq
With $\zeta=1/z$, this means that the zeros of $p$ lie in the set
\begin{eqnarray}
\lcb z \in \complex : \left | \dfrac{1}{z} + \dfrac{\aone}{\azero} \right | \leq \dfrac{\aone}{\azero} + \dfrac{1}{\mutwo} \rcb 
& = & \lcb z \in \complex : \left | 1 + \dfrac{\aone}{\azero} z \right | \leq \lb \dfrac{\aone}{\azero} + \dfrac{1}{\mutwo} \rb |z| \rcb \nonumber \\
& & \nonumber \\
& = & \lcb z \in \complex : \left | z + \dfrac{\azero}{\aone} \right | 
                                              \leq \dfrac{\azero}{\aone} \lb \dfrac{\aone}{\azero} + \dfrac{1}{\mutwo} \rb |z| \rcb \; . \label{recipset}
\end{eqnarray}
The inequality defining the set in~(\ref{recipset}) is of the form
$| z + a  | \leq |a| R |z|$, where $a=\azero/\aone$ and $R=\aone/\azero + 1/\mutwo$. Since $R > 1/|a|$, we can apply Lemma~\ref{lemma_recip} in the 
appendix, which shows that the set defined by~(\ref{recipset}) is the closed exterior of a disk with center $a/(|a|^{2}R^{2}-1)$ and radius
$|a|^{2}R/(|a|^{2}R^{2}-1)$. This means that the zeros of $p$ are excluded from the open disk with center
\beq
\label{centereq}
\dfrac{a}{|a|^{2}R^{2}-1} =  \dfrac{\azero/\aone}{(\azero/\aone)^{2}\lb \aone/\azero+1/\mutwo \rb^{2} - 1}   
                                    =  \dfrac{\mutwo^{2}}{\azero/\aone + 2\mutwo} 
\eeq   
and, from~(\ref{centereq}), with radius 
\bdis
\dfrac{|a|^{2}R}{|a|^{2}R^{2}-1} 
= \bar{a}R \lb \dfrac{a}{|a|^{2}R^{2}-1} \rb
= \dfrac{\azero}{\aone} \lb \dfrac{\aone}{\azero} + \dfrac{1}{\mutwo} \rb \lb \dfrac{\mutwo^{2}}{\azero/\aone + 2\mutwo} \rb
= \dfrac{\mutwo \lb \azero/\aone + \mutwo \rb}{\azero/\aone + 2\mutwo} \; . 
\edis
This completes the proof. \qed 

The first part of Theorem~\ref{Theorem_Bakchik3} was also obtained in Corollary~3 of Theorem~3 in~\cite{AZ}, although it is formulated slightly differently there
with a proof that does not rely on Cauchy-like results, making it more involved algebraically and somewhat difficult to extend beyond the present result. 
No exclusion disk was included in~\cite{AZ}.

The sets $D_{1}$ and $D_{2}$ are typical for the kind of results we will obtain in this section, and for subsequent use we define
\beq
\nonumber
D(a,\gamma) = \lcb z \in \complex : \left | z + a \right | \leq |a| + \gamma \rcb \; , 
\nonumber
\Delta(a,\gamma) = \lcb z \in \complex : \left | z - \dfrac{\gamma^{2}}{|a|+2\gamma}  \right | 
< \dfrac{\gamma \lb |a| + \gamma \rb}{|a| + 2\gamma} \rcb \; . 
\eeq
With this notation, $D_{1}$ and $D_{2}$ in the previous theorem are written as $D_{1} = D(\anmo/\an, \muone)$ and $D_{2}=\Delta(\azero/\aone, \mutwo)$.  

We now derive two more theorems based on Theorem~\ref{theorem_bakchik_cauchy}. Both serve to illustrate the common thread in their
proofs, but each also illustrates a different variation on this theme. The first one shows that the center of the inclusion disk can be varied, while the 
second theorem shows that 
even higher order multipliers, if chosen judiciously, can lead to explicitly computable quantities.

In the following theorem, the inclusion region is determined by a disk centered at $-\varepsilon \anmo/\an$, with $0 < \varepsilon \leq 1$.
%
%
\begin{theorem}
\label{Theorem_Bakchik4a}
Let the real polynomial $p(z) = \sum_{j=0}^{n} a_{j} z^{j}$ with $n \geq 3$ have positive coefficients, let $0 < \varepsilon \leq 1$, and define
\begin{eqnarray*}
& & \gamtwo=\dfrac{(1-\varepsilon) \anmo}{\an} \; , \\
& & \gamone=\dfrac{\anmt-\gamtwo\anmo}{\an} \; , \\ 
& & \gamzero = \max \lcb 0 , \dfrac{\gamone\azero}{-\aone} , \dfrac{\gamtwo\azero+\gamone\aone}{-\atwo} , 
\max_{0 \leq j \leq n-3} \dfrac{\aj-\gamtwo\ajpo-\gamone\ajpt}{\ajpthr} \rcb \; , 
\end{eqnarray*}
and
\begin{eqnarray*}
& & \deltwo=\dfrac{(1-\varepsilon) \aone}{\azero} \; , \\
& & \delone=\dfrac{\atwo-\deltwo\aone}{\azero} \; ,  \\ 
& & \delzero = \max \lcb 0 , \dfrac{\delone\an}{-\anmo} , \dfrac{\deltwo\an+\delone\anmo}{-\anmt} , 
\max_{0 \leq j \leq n-3} \dfrac{\ajpthr-\deltwo\ajpt-\delone\ajpo}{\aj} \rcb \; .
\end{eqnarray*}
Denote by $\muone$ the unique positive zero of $z^{3}-\gamtwo z^{2} -\gamone z -\gamzero$ and by $1/\mutwo$ the unique positive zero of 
$z^{3}-\deltwo z^{2} - \delone z -\delzero$.
Then all the zeros of $p$ lie in the closed disk $D(\varepsilon \anmo/\an, \muone)$, but are excluded from the open disk 
$\Delta(\frac{1}{\varepsilon}\azero/\aone, \mutwo)$.  
\end{theorem}
\prf
Consider $q(z)=(z^{3}-\gamtwo z^{2} -\gamone z -\gamzero)p(z)$ as in the proof of Theorem~\ref{Theorem_Bakchik2}. Then one observes from~(\ref{multipliereq3})
that, with $\gamzero$, $\gamone$, and $\gamtwo$ as in the statement of the current theorem, the second-highest coefficient becomes $\varepsilon\anmo > 0$, 
while all the other coefficients are nonpositive. We therefore conclude, as in the proof of Theorem~\ref{Theorem_Bakchik4}, that all the zeros of $p$ lie
in the closed disk $D(\varepsilon \anmo/\an, \muone)$, where $\muone$ is the unique positive zero of 
$z^{3}-\gamtwo z^{2}-\gamone z -\gamzero$. The open exclusion disk $\Delta(\frac{1}{\varepsilon}\azero/\aone, \mutwo)$ follows analogously as before. \qed  

When $\varepsilon \rightarrow 0^{+}$, Theorem~\ref{Theorem_Bakchik4a} reduces to Theorem~\ref{Theorem_Bakchik2} with the second choice of parameters, and 
when $\varepsilon = 1$, we obtain the following corollary.

%
%
\begin{corollary}    
\label{Theorem_Bakchik4}
Let the real polynomial $p(z) = \sum_{j=0}^{n} a_{j} z^{j}$ with $n \geq 3$ have positive coefficients, and define
\bdis
\gamone= \dfrac{\anmt}{\an} \; , \; 
\gamzero = \max \lcb 0 , \max_{0 \leq j \leq n-3} \dfrac{\aj-\gamone\ajpt}{\ajpthr} \rcb \; , 
\edis
and
\bdis
\delone= \dfrac{\atwo}{\azero} \; , \; 
\delzero = \max \lcb 0 , \max_{0 \leq j \leq n-3} \dfrac{\ajpthr-\delone\ajpo}{\aj} \rcb \; . 
\edis
Denote by $\muone$ the unique positive zero of $z^{3}-\gamone z -\gamzero$ and by $1/\mutwo$ the unique positive zero of $z^{3}-\delone z -\delzero$.
Then all the zeros of $p$ lie in the closed disk $D(\anmo/\an, \muone)$, but are excluded from the open disk $\Delta(\azero/\aone, \mutwo)$. \qed 
\end{corollary}      

Note that $\gamzero =0$ when 
$\anmt/\an = \max_{0 \leq j \leq n-2} \aj/\ajpt$, 
in which case Corollary~\ref{Theorem_Bakchik4} produces the same inclusion disk as in Theorem~\ref{Theorem_Bakchik3}. 
An analogous conclusion follows for the exclusion disk.

The results in the following theorem are obtained by using a quartic multiplier of a particular form that makes it easy to compute its positive zero.
%
%
\begin{theorem}
\label{Theorem_Bakchik5}
Let the real polynomial $p(z) = \sum_{j=0}^{n} a_{j} z^{j}$ with $n \geq 4$ have positive coefficients, and define
\bdis
\alpha= \max \lcb \dfrac{\anmt}{\an} , \dfrac{a_{n-3}}{\anmo} \rcb 
,
\beta = \max \lcb 0 , \max_{0 \leq j \leq n-4} \dfrac{\aj-\alpha\ajpt}{a_{j+4}} \rcb  
,
\muone = \lb \frac{1}{2} \lb \alpha + \lb \alpha^{2} + 4\beta \rb^{1/2} \rb \rb^{1/2} ,
\edis
and
\bdis
\gamma = \max \lcb \dfrac{\atwo}{\azero} , \dfrac{a_{3}}{\aone} \rcb 
,
\delta = \max \lcb 0 , \max_{0 \leq j \leq n-4} \dfrac{a_{j+4}-\gamma \ajpt}{\aj} \rcb  
,
\mutwo = \lb \frac{1}{2} \lb \gamma + \lb \gamma^{2} + 4\delta \rb^{1/2} \rb \rb^{-1/2} ,
\edis
Then all the zeros of $p$ lie in the closed disk $D(\anmo/\an, \muone)$, but are excluded from the open disk $\Delta(\azero/\aone, \mutwo)$.  
\end{theorem}
\prf The proof is similar to that of the previous theorems, except that here we will use a quartic multiplier. Consider
$q(z) = (z^{4} - \alpha z^{2} - \beta)p(z)$:
\begin{multline}
\label{multipliereq6}
q(z) = \an z^{n+4} + \anmo z^{n+3} + \lb \anmt - \alpha \an \rb z^{n+2} + \lb a_{n-3} - \alpha \anmo \rb z^{n+1} \\
+ \sum_{j=4}^{n} \lb a_{j-4} - \alpha \ajmt - \beta \aj  \rb z^{j}  \\ 
- \lb \alpha \aone + \beta a_{3} \rb z^{3} - \lb \alpha \azero + \beta \atwo \rb z^{2} - \beta \aone z - \beta \azero \; .
\end{multline}
Reasoning once again as before, we have from equation~(\ref{multipliereq6}) that all the coefficients of $q$, except those of $z^{n+4}$ and $z^{n+3}$, 
will be nonpositive if 
\bdis
\alpha= \max \lcb \dfrac{\anmt}{\an} , \dfrac{a_{n-3}}{\anmo} \rcb 
\;\; \text{and} \;\;
\beta = \max \lcb 0 , \max_{0 \leq j \leq n-4} \dfrac{\aj-\alpha\ajpt}{a_{j+4}} \rcb \; . 
\edis
Then the Cauchy radius of the second kind of $q$ is the unique positive solution of $(z^{4} - \alpha z^{2} - \beta)p(z)~=~0$,
which is the unique positive solution $\muone$ of $z^{4} - \alpha z^{2} - \beta$.
This quartic is a quadratic in $z^{2}$, and 
its positive zero is easily found to be 
$\muone = \lb \frac{1}{2} \lb \alpha + \lb \alpha^{2} + 4\beta \rb^{1/2} \rb \rb^{1/2}$.
We then have from Theorem~\ref{theorem_bakchik_cauchy} with the set $\Gamma_{2}(1)$ 
that all the zeros of $q$, and therefore also all those of $p$, must lie in the closed disk $D(\anmo/\an, \muone)$.  
As is familiar by now, the open zero exclusion disk $\Delta(\azero/\aone, \mutwo)$ follows from applying to $\recipp$ the result that we just found for $p$. \qed

Here, we observe that, if 
\beq
\label{orstatement}
\dfrac{\anmt}{\an} = \max_{0 \leq j \leq n-2} \dfrac{\aj}{\ajpt} \;\;\;\; \text{or} \;\;\;\; \dfrac{\anmth}{\anmo} = \max_{0 \leq j \leq n-2} \dfrac{\aj}{\ajpt} \; ,
\eeq
then $\beta=0$ and the upper bound becomes $\sqrt{\alpha}$. If, given~(\ref{orstatement}), $\anmt/\an \geq \anmth/\anmo$, 
then we obtain the same bound as in Theorem~\ref{Theorem_Bakchik3}
and Corollary~\ref{Theorem_Bakchik4} while, if, given~(\ref{orstatement}), $\anmth/\anmo \geq \anmt/\an$, then we obtain the same bound as 
in Theorem~\ref{Theorem_Bakchik3}, but not necessarily as in Corollary~\ref{Theorem_Bakchik4}. Once again, analogous statements hold for the exclusion disk.

\begin{figure}[H]
\begin{center}
\raisebox{0ex}{\includegraphics[width=0.65\linewidth]{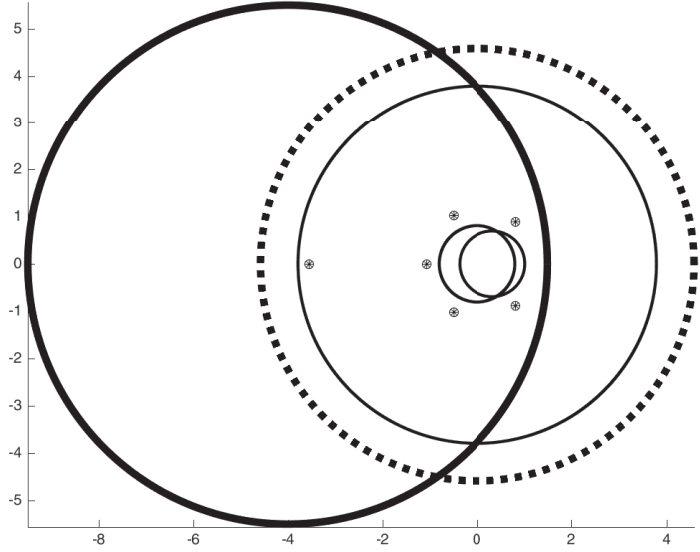}}
\caption{Inclusion regions from Theorems~\ref{Theorem_Bakchik2}(2) and \ref{Theorem_Bakchik5} for $p_{3}$.}
\label{eek_fig1}                      
\end{center}
\end{figure}

Let us now graphically illustrate some of these results.
\vskip 0.25 cm \noindent {\bf Example.}
Since it would be impractical to form all possible combinations, we have chosen the upper and lower bounds bounds from Theorem~\ref{Theorem_Bakchik2}(2)  
(with the second choice of parameters) and combined them with the inclusion and exclusion disks from Theorem~\ref{Theorem_Bakchik5} for the polynomial
$p_{3}(z)=z^{6}+4z^{5}+2z^{4}+2z^{3}+3z^{2}+6z+7$ we encountered at the end of Section~\ref{onedisk_origin}.
The result can be found in Figure~\ref{eek_fig1}, where the solid circles centered at the 
origin (right-most circles) are obtained from Theorem~\ref{Theorem_Bakchik2}(2), 
and the others from Theorem~\ref{Theorem_Bakchik5}. The dashed circle indicates the upper bound
from the Cauchy radius. The circle obtained from the Cauchy radius of the second kind almost coincides with the one from Theorem~\ref{Theorem_Bakchik5},
and is therefore not drawn.
The zeros of $p$ are indicated by circled asterisks. Here the disk from Theorem~\ref{Theorem_Bakchik5} cuts off a significant part of the disk from 
Theorem~\ref{Theorem_Bakchik2}.

The radii of the exclusion and inclusion disks, respectively, for the various theorems used in Figure~\ref{eek_fig1} are as follows:
\bdis
\begin{array}{ll}
\text{Cauchy radii of the second kind:}  &   0.698 , 5.475  \; , \\
\text{Theorem~\ref{Theorem_Bakchik3}:}  & 0.513 , 5.732  \; ,  \\
\text{Corollary~\ref{Theorem_Bakchik4}:}  &  0.607 ,  5.618\; ,  \\
\text{Theorem~\ref{Theorem_Bakchik5}:}  & 0.693 , 5.492 \; .  \\
\end{array}
\edis
All inclusion disks have the same center ($-\anmo/\an$), but the exclusion disks do not. The corresponding radii for the disks from 
Theorem~\ref{Theorem_Bakchik4a}, where the inclusion disk is centered at $-\frac{1}{2}\anmo/an$ (obtained with $\varepsilon=1/2$), 
are $0.664$ and $4.740$. As in Section~\ref{onedisk_origin},    
one discerns a positive trend with increasing degree of the multipliers upon which the theorems are based. 

%
%

We now numerically compare the radii of the disks obtained in the theorems of this section using the same two classes of polynomials 
as in Section~\ref{onedisk_origin}. Here too, we generated $1000$ polynomials for each class and collected the results in the following two tables.
\begin{itemize}
\item Table~\ref{Table_21} lists, for Class I polynomials, the median of the radii of the disks centered at $-\anmo/\an$.
\item Table~\ref{Table_22} is the analog of Table~\ref{Table_21} for Class II polynomials. 
\end{itemize}
We have designated by "Cauchy" the radius obtained from the set $\Gamma_{2}(1)$ in Theorem~\ref{theorem_bakchik_cauchy}, i.e., $s_{2}+\anmo/\an$,
where $s_{2}$ is the Cauchy radius of the second kind.

%
%
\begin{table}[H]
\begin{center}
\small           
\begin{tabular}{c|c|c|c|c}
                  & Cauchy               & Theorem 4.1          &  Corollary 4.1       & Theorem 4.3         \\ \hline         
                  &                      &                      &                      &                     \\
n=10              &        2.457         &        3.957         &        2.808         &       2.741         \\        
                  &                      &                      &                      &                     \\  
n=40              &        2.471         &        7.235         &        3.884         &       3.447         \\        
\end{tabular}
\caption{Comparison of inclusion regions from Section~\ref{onedisk_notorigin} for Class I polynomials.}
\label{Table_21}
\end{center}
\end{table}
\normalsize

%
%
\begin{table}[H]
\begin{center}
\small           
\begin{tabular}{c|c|c|c|c}
                  & Cauchy               & Theorem 4.1          &  Corollary 4.1       & Theorem 4.3         \\ \hline         
                  &                      &                      &                      &                     \\
n=10              &        2.413         &        2.649         &        2.362         &       2.359         \\        
                  &                      &                      &                      &                     \\
n=40              &        2.403         &        2.857         &        2.463         &       2.411         \\        
\end{tabular}
\caption{Comparison of inclusion regions from Section~\ref{onedisk_notorigin} for Class II polynomials.}
\label{Table_22}
\end{center}
\end{table}
\normalsize

We observe that, for both classes of polynomials, Corollary~\ref{Theorem_Bakchik4} and Theorem~\ref{Theorem_Bakchik5} outperform Theorem~\ref{Theorem_Bakchik3}.
Corollary~\ref{Theorem_Bakchik4} is better than Theorem~\ref{Theorem_Bakchik5} for Class~I polynomials, but the reverse is true for Class~II polynomials 
with increasing degree of the polynomials. The disk based on the Cauchy radius of the second kind (which requires the solution of a polynomial equation) 
is generally smaller than any of the disks obtained here, especially when the degree of the polynomial increases.

So far, we have only used the sets $\Gamma_{1}(k)$ and $\Gamma_{2}(k)$ from Theorem~\ref{theorem_bakchik_cauchy} for $k=1$. In the following two sections we 
derive sets for $k \geq 2$ consisting of two and three disks, respectively. We limit ourselves to zero inclusion regions obtained using the polynomial, but 
not the reverse polynomial. 

%
%
\section{Two disks}           
\label{twodisks}         

In this section we derive zero inclusion regions consisting of two disks, organized in three theorems. The first of these relies on a quadratic multiplier,
while the other two relie on cubic multipliers.

%
%
\begin{theorem}
\label{theorem_lemniscates1}          
Let the real polynomial $p(z) = \sum_{j=0}^{n} a_{j} z^{j}$ with $n \geq 2$ have positive coefficients, and define
\bdis
\mu = \lb \max_{0 \leq j \leq n-2} \dfrac{\aj}{\ajpt} \rb^{1/2}         
\text{and} \;\;
R= \mu^{2}+\dfrac{\anmo}{\an} \mu \; .
\edis

\noindent Then all the zeros of $p$ are included in the union of disks 
\beq
\nonumber
\lcb z \in \complex : |z| \leq R^{1/2} \rcb 
\bigcup 
\lcb z \in \complex : \left | z + \dfrac{\anmo}{\an} \right | \leq R^{1/2} \rcb \; . 
\eeq
If the disks are disjoint, then 
the disk centered at $-\anmo/\an$ contains one zero and the one centered at the origin contains the remaining $n-1$ zeros of $p$.
\end{theorem}
\prf
Consider, as in the proof of Theorem~\ref{Theorem_Bakchik3}, the polynomial $q(z)=(z^{2}-\gamma)p(z)$. Then 
\bdis
q(z) = \an z^{n+2} + \anmo z^{n+1} + \sum_{j=2}^{n} (\ajmt - \gamma \aj ) z^{j} -\gamma \aone z -\gamma \azero \; . 
\edis
If we now choose $\gamma=\mu^{2}$ with $\mu$ as in the statement of the theorem, then all the coefficients of $q$, other than the two leading ones,
are nonpositive, so that $s_{2}$, its Cauchy radius of the second kind is its unique positive solution, namely, $\mu$. The set $\Gamma_{1}(2)$ in
Theorem~\ref{theorem_bakchik_cauchy} is then given by 
\bdis
\Gamma_{1}(2) = \lcb z \in \complex : \left | z \lb z + \dfrac{\anmo}{\an} \rb \right | \leq R \rcb \; , 
\edis
with $R=\mu^{2}+(\anmo/\an)\mu$. This set is contained in the following union of two disks: 
\bdis
\lcb z \in \complex : |z| \leq R^{1/2} \rcb 
\bigcup 
\lcb z \in \complex : \left | z + \dfrac{\anmo}{\an} \right | \leq R^{1/2} \rcb \; . 
\edis

If the disks are disjoint, then from Theorem~\ref{theorem_bakchik_cauchy} with $k=2$, we obtain that
the disk centered at the origin contains $1+(n+2)-2=n+1$ zeros, while the disk centered at $-\anmo/\an$ contains one zero of $q$. 
The zeros of $q$ are those of $p$ with the addition of $\pm \mu$. Because the disks are disjoint, all points in the disk centered at $-\anmo/\an$ have a 
negative real part, so that $\mu$ must lie in the disk centered at the origin, which means that $-\mu$ also lies in that disk. The only zero of $q$
in the disk centered at $-\anmo/\an$ must therefore be a zero of $p$, while the other $n-1$ zeros of $p$ lie in the disk centered at the origin.
This concludes the proof. \qed

This theorem is closely related to Theorem~\ref{Theorem_Bakchik3}. However, instead of one disk centered at $-\anmo/\an$, we now have two 
disks - one centered at the origin and the other at $-\anmo/\an$ - with a smaller radius. That this radius is smaller follows from the fact that 
\bdis
\mu^{2}+\dfrac{\anmo}{\an} \mu \leq \lb \mu+\dfrac{\anmo}{\an} \rb^{2} 
\Longrightarrow
\lb \mu^{2}+\dfrac{\anmo}{\an} \mu \rb^{1/2} \leq \mu+\dfrac{\anmo}{\an} \; ,
\edis
where $\mu$ has the same meaning here as in Theorem~\ref{Theorem_Bakchik3}.

%
%
\begin{theorem}
\label{theorem_lemniscates2}          
Let the real polynomial $p(z) = \sum_{j=0}^{n} a_{j} z^{j}$ with $n \geq 3$ have positive coefficients, let
$0 < \varepsilon \leq 1$, and define 
\bdis
\gamma_{2} = \dfrac{(1-\varepsilon)\anmo}{\an} \; , \; \gamma_{1} = \dfrac{\anmt - \gamma_{2} \anmo}{\an} \; , 
\edis
and
\bdis
\gamma_{0} =  \max \lcb 0 , \dfrac{\gamone\azero}{-\aone} , \dfrac{\gamtwo\azero+\gamone\aone}{-\atwo} , 
\max_{0 \leq j \leq n-3} \dfrac{\aj-\gamma_{2}\ajpo-\gamma_{1}\ajpt}{a_{j+3}} \rcb \; .
\edis
Let $\mu$ be the unique positive zero of $z^{3}-\gamma_{2}z^{2}-\gamma_{1}z-\gamma_{0}$, and let $R=\mu^{3}+(\varepsilon\anmo/\an)\mu^{2}$.  

\noindent Then all the zeros of $p$ are included in the union of disks 
\beq
\nonumber
\lcb z \in \complex : |z| \leq R^{1/3} \rcb 
\bigcup 
\lcb z \in \complex : \left | z + \varepsilon \dfrac{\anmo}{\an} \right | \leq R^{1/3} \rcb \; . 
\eeq
If the disks are disjoint, then the disk centered at $-\varepsilon \anmo/\an$ contains one zero and the one centered at the origin contains the 
remaining $n-1$ zeros of $p$.
\end{theorem}
\prf
As in the proof of Theorems~\ref{Theorem_Bakchik2} and~\ref{Theorem_Bakchik4a}, consider the polynomial $q(z)=(z^{3}-\gamtwo z^{2}-\gamone z -\gamzero)p(z)$:
\begin{multline}
\nonumber                          
q(z) = \an z^{n+3} + (\anmo -\gamtwo \an) z^{n+2} + (\anmt - \gamtwo \anmo - \gamone \an )z^{n+1} \\ 
\hskip 2cm + \sum_{j=3}^{n} \lb \ajmthr - \gamtwo \ajmt - \gamone \ajmo - \gamzero \aj \rb  z^{j} 
- \lb \gamtwo \azero +\gamone \aone +\gamzero\atwo \rb z^{2}  \\
\hskip 3cm - \lb \gamone \azero +\gamzero \aone \rb z -\gamzero\azero \; . 
\end{multline}
If we choose $\gamma_{0}$, $\gamma_{1}$, and $\gamma_{2}$ as in the statement of the theorem, then the two leading coefficents are positive, the coefficient
of $z^{n+1}$ vanishes and all other coefficients are nonpositive. The Cauchy radius of the third kind of $q$ is its unique positive zero, which is the unique
positive zero $\mu$ of $z^{3}-\gamtwo z^{2}-\gamone z -\gamzero$. This means that the set $\Gamma_{1}(3)$ in Theorem~\ref{theorem_bakchik_cauchy}
is given by 
\bdis
\Gamma_{1}(3) = \lcb z \in \complex : \left | z^{2} \lb z + \varepsilon \dfrac{\anmo}{\an} \rb \right | \leq R \rcb \; , 
\edis
where $R=\mu^{3}+(\anmo/\an)\mu^{2}$. This set is contained in the following union of two disks: 
\bdis
\lcb z \in \complex : |z| \leq R^{1/3} \rcb 
\bigcup 
\lcb z \in \complex : \left | z + \varepsilon \dfrac{\anmo}{\an} \right | \leq R^{1/3} \rcb \; . 
\edis

If the disks are disjoint, then from Theorem~\ref{theorem_bakchik_cauchy} with $k=3$, we have that
the disk centered at the origin contains $2+(n+3)-3=n+2$ zeros, while the disk centered at $-\varepsilon\anmo/\an$ contains one zero of $q$. 
The zeros of $q$ are those of $p$ with the addition of the three zeros of $h(z) = z^{3}-\gamtwo z^{2}-\gamone z -\gamzero$. This polynomial $h$ has a unique 
positive zero, which we denoted by $\mu$, and its coefficients $\gamma_{2}$ and $\gamma_{0}$ are nonnegative. Let us first assume that $\gamma_{1}$ is also 
nonnegative. Then $\mu$ has the largest modulus of all the zeros of $h$. Because the disks are disjoint, all points in the disk centered at 
$-\varepsilon \anmo/\an$ have 
a negative real part, so that $\mu$ must lie in the disk centered at the origin. Since it has the largest modulus of all zeros of $h$, the other two zeros must 
also lie in that disk, so that the only zero of $q$ in the disk centered at $-\varepsilon \anmo/\an$ must be a zero of $p$, 
while the other $n-1$ zeros of $p$ lie in the disk 
centered at the origin.

If $\gamma_{1} < 0$, then $h$ has no negative zeros, since these negative zeros are the positive zeros of 
\bdis
(-\zeta)^{3} -\gamma_{2}(-\zeta)^{2}-\gamma_{1}(-\zeta)-\gamma_{0}
=
-\zeta^{3} -\gamma_{2}\zeta^{2}+\gamma_{1}\zeta-\gamma_{0} \; ,
\edis
and this is a polynomial with nonpositive coefficients, which does not have positive zeros. The other two zeros of $h$ must therefore be complex
(and each other's conjugate).
If we denote these complex zeros by $\nu$ and $\bar{\nu}$, then we have that 
\beq
\label{munuineq}
\mu + \nu + \bar{\nu} = \gamma_{2} \Longrightarrow \mu + 2 \text{Re}(\nu) > 0 \; .
\eeq
If $\text{Re}(\nu) \geq 0$, then $\nu$ and $\bar{\nu}$ cannot lie in the disk centered at $-\varepsilon \anmo/\an$, and must therefore be in the disk centered 
at the origin. If $\text{Re}(\nu) < 0$, then from~(\ref{munuineq}), we have that $\text{Re}(\nu) > -\mu/2$, which means that $\nu$ and $\bar{\nu}$ must lie in the
disk centered at the origin because that disk contains $\mu$. 
We conclude that the only zero of $q$ in the disk centered at $-\varepsilon \anmo/\an$ must be a zero of $p$. The remaining
$n-1$ zeros of $p$ are contained in the disk centered at the origin. This concludes the proof. \qed

When $\varepsilon \rightarrow 0^{+}$, Theorem~\ref{theorem_lemniscates2} reduces to Theorem~\ref{Theorem_Bakchik2} with the second choice of parameters, and 
when $\varepsilon = 1$, we obtain the following corollary.

%
%
\begin{corollary}
\label{corollary_lemniscates2}          
Let the real polynomial $p(z) = \sum_{j=0}^{n} a_{j} z^{j}$ with $n \geq 3$ have positive coefficients.
Define 
\bdis
\gamma_{1} = \dfrac{\anmt}{\an} 
\;\; \text{and} \;\;
\gamma_{0} =  \max \lcb 0 , \max_{0 \leq j \leq n-3} \dfrac{\aj-\gamma_{1}\ajpt}{a_{j+3}} \rcb \; .
\edis
Let $\mu$ be the unique positive zero of $z^{3}-\gamma_{1}z-\gamma_{0}$, and let $R=\mu^{3}+(\anmo/\an)\mu^{2}$.  

\noindent Then all the zeros of $p$ are included in the union of disks 
\beq
\nonumber
\lcb z \in \complex : |z| \leq R^{1/3} \rcb 
\bigcup 
\lcb z \in \complex : \left | z + \dfrac{\anmo}{\an} \right | \leq R^{1/3} \rcb \; . 
\eeq
If the disks are disjoint, then
the disk centered at $-\anmo/\an$ contains one zero and the one centered at the origin contains the remaining $n-1$ zeros of $p$.
\end{corollary}     

\vskip 0.25cm {\bf Remarks.}
Theorem~\ref{theorem_lemniscates2} is similar to Theorem~\ref{Theorem_Bakchik4a}. However, we now have two disks instead of one, both
with the same radius, and this radius is smaller than that of the disk in Theorem~\ref{Theorem_Bakchik4a} because                   
\beq
\label{smallerradius}
\mu^{3} + \varepsilon \dfrac{\anmo}{\an} \mu^{2}  \leq \lb \mu + \varepsilon \dfrac{\anmo}{\an} \rb^{3} 
\Longrightarrow
\lb \mu^{3} + \varepsilon \dfrac{\anmo}{\an} \mu^{2} \rb^{1/3} \leq \mu + \varepsilon \dfrac{\anmo}{\an} \; .
\eeq

A similar observation holds for Corollary~\ref{corollary_lemniscates2} and Corollary~\ref{Theorem_Bakchik4}.

%
%
\begin{theorem}
\label{theorem_lemniscates3}          
Let the real polynomial $p(z) = \sum_{j=0}^{n} a_{j} z^{j}$ with $n \geq 3$ have positive coefficients, and define
\bdis
\mu = \lb \max_{0 \leq j \leq n-3} \dfrac{\aj}{a_{j+3}} \rb^{1/3}         
\text{and} \;\;
R= \mu^{2}+\dfrac{\anmo}{\an} \mu + \dfrac{\anmt}{\an} \; ,
\edis
and let $c_{1}$ and $c_{2}$ be the zeros of the quadratic $z^{2}+(\anmo/\an) z + \anmt/\an$.

\noindent Then all the zeros of $p$ are included in the union of disks 
\beq
\nonumber
\lcb z \in \complex : \left |z - c_{1} \right | \leq R^{1/2} \rcb 
\bigcup 
\lcb z \in \complex : \left | z -c_{2} \right | \leq R^{1/2} \rcb \; . 
\eeq

If the disks are disjoint, then 
$c_{2} < c_{1} < 0$, the disk centered at $c_{2}$ contains one zero, and the one centered at $c_{1}$ contains the remaining $n-1$ zeros of $p$.
\end{theorem}
\prf
Let the polynomial $q$ be defined by $q(z)=(z^{3}-\gamma)p(z)$. Then 
\beq                  
\nonumber                          
q(z) = \an z^{n+3} + \anmo z^{n+2} + \anmt z^{n+1} + \sum_{j=3}^{n} \lb \ajmthr - \gamma \aj \rb  z^{j} 
         - \gamma\atwo  z^{2} - \gamma \aone z -\gamma\azero \; . 
\eeq
If $\gamma=\mu^{3}$ with $\mu$ as in the statement of the theorem, then all the coefficients of $q$, other than the three leading ones,
are nonpositive, so that its Cauchy radius of the third kind is its unique positive zero, namely, $\mu$. The set $\Gamma_{2}(2)$ in
Theorem~\ref{theorem_bakchik_cauchy} is then given by 
\bdis
\Gamma_{2}(2) = \lcb z \in \complex : \left | z^{2} + \dfrac{\anmo}{\an} z + \dfrac{\anmt}{\an} \right | \leq R \rcb \; , 
\edis
with $R=\mu^{2}+(\anmo/\an)\mu+\anmt/\an$. Let 
\bdis
c_{1} = \dfrac{1}{2} \lb -\dfrac{\anmo}{\an} + \lb \dfrac{\anmo^{2}}{\an^{2}} - 4 \dfrac{\anmt}{\an} \rb ^{1/2} \rb
\;\; \text{and} \;\; 
c_{2} = \dfrac{1}{2} \lb -\dfrac{\anmo}{\an} - \lb \dfrac{\anmo^{2}}{\an^{2}} - 4 \dfrac{\anmt}{\an} \rb ^{1/2} \rb
\edis
be the zeros of $z^{2} + (\anmo/\an) z + \anmt/\an$, then this set is contained in the following union of two disks: 
\beq
\nonumber
\lcb z \in \complex : \left |z - c_{1} \right | \leq R^{1/2} \rcb 
\bigcup 
\lcb z \in \complex : \left | z -c_{2} \right | \leq R^{1/2} \rcb \; . 
\eeq

These disks are disjoint if and only if 
\beq    
\nonumber
\left | c_{1}-c_{2} \right | = \left | \lb \dfrac{\anmo^{2}}{\an^{2}} - 4 \dfrac{\anmt}{\an} \rb^{1/2} \right | > 2 R^{1/2} \; .      
\eeq
If $\anmo^{2}/\an^{2} \leq 4 \anmt/\an$, then 
\bdis
\left | c_{1}-c_{2} \right | = \left | \lb \dfrac{\anmo^{2}}{\an^{2}} - 4 \dfrac{\anmt}{\an} \rb^{1/2} \right |      
= \left | \lb 4 \dfrac{\anmt}{\an} - \dfrac{\anmo^{2}}{\an^{2}} \rb^{1/2} i \right |      
= \lb 4 \dfrac{\anmt}{\an} - \dfrac{\anmo^{2}}{\an^{2}} \rb^{1/2} \; .  
\edis
In this case, it is impossible that $\left | c_{1} - c_{2} \right | > 2R^{1/2}$ since 
$\left | c_{1} - c_{2} \right | < 2(\anmt/\an)^{1/2}$, while $2R^{1/2} > 2(\anmt/\an)^{1/2}$. 
This means that when the disks are disjoint, then $\anmo^{2}/\an^{2} > 4 \anmt/\an$ and both $c_{1}$ and $c_{2}$
are real and negative with $c_{2} < c_{1} < 0$. 

The disk centered at $c_{1}$ must contain the origin since otherwise, by Theorem~\ref{theorem_bakchik_cauchy}, 
the set $\Gamma_{2}(2)$ would only contain two zeros of $q$, when it has, in fact, $n+3$ zeros. This can also be shown explicitly.
As a consequence, the disk centered at $c_{1}$ contains $1+(n+3)-2=n+2$ zeros, while the disk centered at $c_{2}$ contains one zero of $q$.
Since $\mu > 0$, it must lie in the disk centered at $c_{1}$ and because the other two zeros of $z^{3}-\mu^{3}$ both have the same modulus $\mu$, they
too must lie in this disk. Therefore, the single zero of $q$ that lies in the disk centered at $c_{2}$ cannot be a zero of $z^{3}-\mu^{3}$ and must
therefore be a zero of $p$. Clearly, the remaining $n-1$ zeros of $p$ lie in the disk centered at $c_{1}$.
This concludes the proof. \qed

It is generally difficult to predict which theorem produces the better result, but to obtain some idea, we numerically compare the radii of the disks 
from Theorem~\ref{theorem_lemniscates1} and Corollary~\ref{corollary_lemniscates2}, as they have the same centers. 
We generated $1000$ polynomials for each of the same two classes as in Section~\ref{onedisk_origin}, and listed the results 
in the following two tables.
\begin{itemize}
\item Table~\ref{Table_31} lists, for Class I polynomials, the median of the radii of the disks centered at the origin and $-\anmo/\an$. 
\item Table~\ref{Table_32} is the analog of Table~\ref{Table_31} for Class II polynomials. 
\end{itemize}
We have designated by "Cauchy" the radius obtained from the set $\Gamma_{1}(2)$ in Theorem~\ref{theorem_bakchik_cauchy}, i.e., 
$\lb s_{2}^{2}+(\anmo/\an) s_{2} \rb^{1/2}$, where $s_{2}$ is the Cauchy radius of the second kind.

%
%
\begin{table}[H]
\begin{center}
\small           
\begin{tabular}{c|c|c|c|c}
                  &       Cauchy         &     Theorem 5.1         & Corollary 5.1        \\ \hline         
                  &                      &                         &                      \\
n=10              &        1.845         &        3.315            &    2.052             \\        
                  &                      &                         &                      \\  
n=40              &        1.846         &        6.292            &    2.884             \\        
\end{tabular}
\caption{Comparison of inclusion regions from Section~\ref{twodisks} for Class I polynomials.}
\label{Table_31}
\end{center}
\end{table}
\normalsize

%
%
\begin{table}[H]
\begin{center}
\small           
\begin{tabular}{c|c|c|c|c}
                  &        Cauchy        &    Theorem 5.1       &  Corollary 5.1     \\ \hline         
                  &                      &                      &                      \\
n=10              &        1.845         &        2.072         &      1.633           \\        
                  &                      &                      &                      \\
n=40              &        1.854         &        2.360         &      1.795           \\        
\end{tabular}
\caption{Comparison of inclusion regions from Section~\ref{twodisks} for Class II polynomials.}
\label{Table_32}
\end{center}
\end{table}
\normalsize

\begin{figure}[H]
\begin{center}
\raisebox{0ex}{\includegraphics[width=0.475\linewidth]{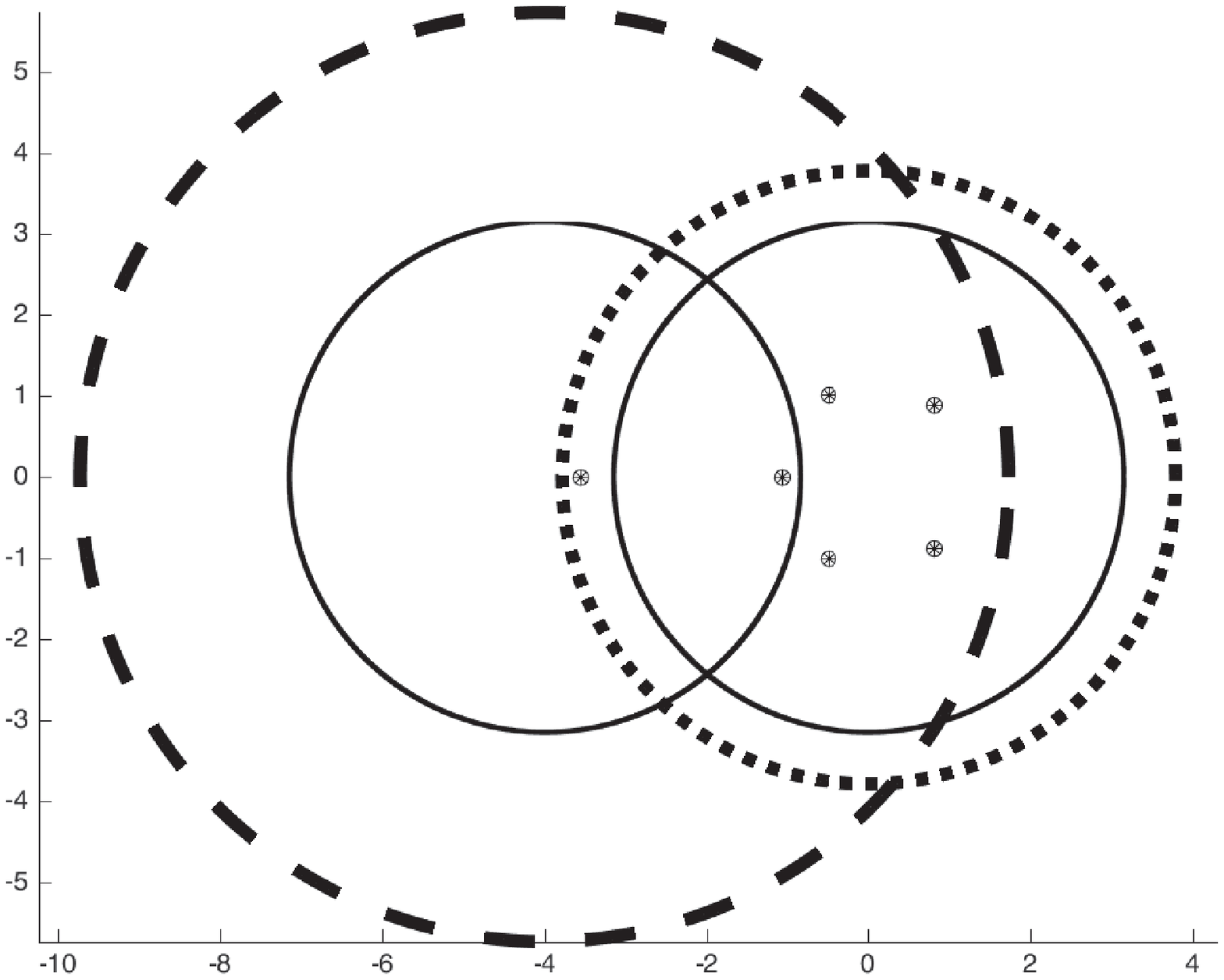}}
\caption{Comparison of Theorems~\ref{Theorem_Bakchik3} and \ref{theorem_lemniscates1} for $p_{3}$.}
\label{eek_fig2}                      
\end{center}
\end{figure}

We observe that, for both classes of polynomials, Corollary~\ref{corollary_lemniscates2} generally outperforms Theorem~\ref{theorem_lemniscates1}, 
although on rare occasions the opposite is true.
The disks based on the Cauchy radius of the second kind, which, we recall, requires the solution of a polynomial equation, are generally 
smaller than any of the disks obtained here for Class I polynomials, but not for Class II polynomials, where Corollary~\ref{corollary_lemniscates2}
clearly delivers better results on average.
\vskip 0.25cm
\noindent {\bf Example.} 
Figure~\ref{eek_fig2} compares the inclusion regions obtained from Theorem~\ref{Theorem_Bakchik3} and Theorem~\ref{theorem_lemniscates1} for the same 
polynomial $p_{3}(z)=z^{6}+4z^{5}+2z^{4}+2z^{3}+3z^{2}+6z+7$ that we used at the end of Section~\ref{onedisk_origin} and in 
Figure~\ref{eek_fig1}.
The dotted circle centered at the origin (with radius 3.788) is the boundary of the disk obtained from Theorem~\ref{Theorem_Bakchik2}(2), while the large 
dashed circle (with radius 5.732) is the
boundary of the disk from Theorem~\ref{Theorem_Bakchik3}. The two smaller solid circles (with radii 3.151) are the boundaries of the disks from 
Theorem~\ref{theorem_lemniscates1}.
The zeros of $p$ are indicated by circled asterisks. The Cauchy radius of $p_{3}$ is 4.580. 

%
%

Regions composed of two disks are not necessarily smaller than those composed of a single disk, although they frequently are.
Moreover, when the disks are disjoint, they provide additional information about the location of the zeros that cannot be obtained from 
standard generalizations of the \eneskakk theorem.

\section{Three disks}           
\label{threedisks}         

In this section we carry out one more application of Theorem~\ref{theorem_bakchik_cauchy} to obtain zero inclusion regions consisting
of three disks.

%
%
\begin{theorem}
\label{theorem_lemniscates4}          
Let the real polynomial $p(z) = \sum_{j=0}^{n} a_{j} z^{j}$ with $n \geq 3$ have positive coefficients.
Let $0 < \varepsilon \leq 1$, and define 
\bdis
\gamma_{1} = \dfrac{(1-\varepsilon)\anmt}{\an} 
\;\; \text{and} \;\;
\gamma_{0} =  \max \lcb 0 , \max_{0 \leq j \leq n-3} \dfrac{\aj-\gamma_{1}\ajpt}{a_{j+3}} \rcb \; .
\edis
Let $\mu$ be the unique positive zero of $z^{3}-\gamma_{1}z-\gamma_{0}$, let $R=\mu^{3}+(\anmo/\an)\mu^{2}+\varepsilon(\anmt/\an)\mu$,  
and let $c_{1}(\varepsilon)$ and $c_{2}(\varepsilon)$ be the zeros of the quadratic $z^{2}+(\anmo/\an) z + \varepsilon\anmt/\an$.

\noindent Then all the zeros of $p$ are included in the union of disks 
\beq
\nonumber
\lcb z \in \complex : |z| \leq R^{1/3} \rcb
\bigcup 
\lcb z \in \complex : \left | z - c_{1} (\varepsilon) \right | \leq R^{1/3} \rcb 
\bigcup 
\lcb z \in \complex : \left | z - c_{2} (\varepsilon) \right | \leq R^{1/3} \rcb \; . 
\eeq
There exist only the following two scenarios for disks to be disjoint.
\begin{enumerate}
\item[(1)] 
The disk centered at the origin is disjoint from the other two, in which case that disk contains 
$n-2$ zeros of $p$, while the union of the other two contains the two remaining zeros of $p$. If these are also disjoint, then each contains one zero of $p$. 
\item[(2)] 
The two disks not centered at the orgin are disjoint, but  
only one of them is disjoint from the disk at the origin, in which case that disk contains one zero of $p$, while the union of the other two contains 
$n-1$ zeros of $p$. This scenario is only possible when $c_{1}(\varepsilon)$ and $c_{2}(\varepsilon)$ are real and negative.
\end{enumerate}
\end{theorem}
\prf
Consider the polynomial $q(z)=z^{3}-\gamma_{1} z-\gamma_{0}$, given by
\begin{multline}
\nonumber                          
q(z) = \an z^{n+3} + \anmo z^{n+2} + (\anmt - \gamone \an )z^{n+1} \\ 
\hskip 2cm + \sum_{j=3}^{n} \lb \ajmthr - \gamone \ajmo - \gamzero \aj \rb  z^{j} 
- \lb \gamone \aone +\gamzero\atwo \rb z^{2}  \\
\hskip 3cm - \lb \gamone \azero +\gamzero \aone \rb z -\gamzero\azero \; . 
\end{multline}
If we choose $\gamma_{0}$ and $\gamma_{1}$ as in the statement of the theorem, then the three leading coefficents are positive, while all other
coefficients are nonpositive. The Cauchy radius of the third kind of $q$ is then its unique positive zero, which is the unique
positive zero $\mu$ of $z^{3}-\gamone z -\gamzero$. This means that the set $\Gamma_{1}(3)$ in Theorem~\ref{theorem_bakchik_cauchy}
is given by 
\bdis
\Gamma_{1}(3) = \lcb z \in \complex : \left | z \lb z^{2} + \dfrac{\anmo}{\an} z + \varepsilon \dfrac{\anmt}{\an}  \rb \right | \leq R \rcb \; , 
\edis
where $R=\mu^{3}+(\anmo/\an)\mu^{2}+\varepsilon (\anmt/\an)\mu$. This set is contained in the following union of three disks: 
\bdis
\lcb z \in \complex : |z| \leq R^{1/3} \rcb
\bigcup 
\lcb z \in \complex : \left | z - c_{1} (\varepsilon) \right | \leq R^{1/3} \rcb 
\bigcup 
\lcb z \in \complex : \left | z - c_{2} (\varepsilon) \right | \leq R^{1/3} \rcb \; . 
\edis

Several scenarios arise when the disks are disjoint. 
If 
\bdis
\dfrac{\anmo^{2}}{\an^{2}} \leq 4 \dfrac{\anmt}{\an} \; ,
\edis
then $c_{1}(\varepsilon)$ and $c_{1}(\varepsilon)$ are complex conjugate with a negative real part, and the disks centered at these points are 
either not disjoint or both disjoint from the disk centered at the origin. When they are disjoint, then $q$ has $1+(n+3)-3=n+1$ zeros in the disk 
centered at the orgin. One of those must be $\mu$,
which is the real positive zero of $h(z):=z^{3}-\gamone z - \gamzero$ with largest modulus. The other two zeros of $h$ must therefore also lie in that disk, 
and this means that the remaining $n-2$ zeros of $q$ in that disk are zeros of $p$, while two zeros of $p$ lie in the union of the two disks centered 
at $c_{1}(\varepsilon)$ and $c_{2}(\varepsilon)$. If these are also disjoint from each other, then each contains one zero of $p$. 

If $c_{1}(\varepsilon)$ and $c_{1}(\varepsilon)$ are not complex, then they are both real and negative. If the disks centered at these points are
both disjoint from the disk centered at the origin, then reasoning similarly as before, their union contains two zeros of $p$; if they are
disjoint from each other, then each contains one zero of $p$. If only one is disjoint from the disk centered at the origin, then it contains one zero
of $p$. \qed

When $\varepsilon \rightarrow 0^{+}$, Theorem~\ref{theorem_lemniscates4} reduces to Corollary~\ref{corollary_lemniscates2}, and 
when $\varepsilon = 1$, we obtain the following corollary.

%
%
\begin{corollary}  
\label{corollary_lemniscates4}          
Let the real polynomial $p(z) = \sum_{j=0}^{n} a_{j} z^{j}$ with $n \geq 3$ have positive coefficients.
Define 
\bdis
\mu = \lb \max_{0 \leq j \leq n-3} \dfrac{\aj}{a_{j+3}} \rb^{1/3}         
\text{and} \;\;
R= \mu^{3}+\dfrac{\anmo}{\an} \mu^{2} + \dfrac{\anmt}{\an} \mu \; ,
\edis
and let $c_{1}$ and $c_{2}$ be the zeros of the quadratic $z^{2}+(\anmo/\an) z + \anmt/\an$.

\noindent Then all the zeros of $p$ are included in the union of disks 
\beq
\nonumber
\lcb z \in \complex : |z| \leq R^{1/3} \rcb
\bigcup 
\lcb z \in \complex : \left |z - c_{1} \right | \leq R^{1/3} \rcb 
\bigcup 
\lcb z \in \complex : \left | z -c_{2} \right | \leq R^{1/3} \rcb \; . 
\eeq
There exist only the following two scenarios for disks to be disjoint.
\begin{enumerate}
\item[(1)] 
The disk centered at the origin is disjoint from the other two, in which case that disk contains 
$n-2$ zeros of $p$, while the union of the other two contains the two remaining zeros of $p$. If these are also disjoint, then each contains one zero of $p$. 
\item[(2)] 
The two disks not centered at the orgin are disjoint, but  
only one of them is disjoint from the disk at the origin, in which case that disk contains one zero of $p$, while the union of the other two contains 
$n-1$ zeros of $p$. This scenario is only possible when $c_{1}$ and $c_{2}$ are real and negative.
\end{enumerate}
\end{corollary}            

\noindent {\bf Remarks.}
\begin{itemize}
\item
We remark that the two disks in Corollary~\ref{corollary_lemniscates4} that are not centered at the origin have the same centers as 
the disks in Theorem~\ref{theorem_lemniscates3}, but their radii are smaller. This is so because
\bdis
\mu^{2} \lb \mu^{2}+\dfrac{\anmo}{\an} \mu + \dfrac{\anmt}{\an} \rb^{2} \leq \lb \mu^{2}+\dfrac{\anmo}{\an} \mu + \dfrac{\anmt}{\an} \rb^{3}  \\
\edis
implies that
\bdis
\mu^{1/3} \lb \mu^{2}+\dfrac{\anmo}{\an} \mu + \dfrac{\anmt}{\an} \rb^{1/3} \leq \lb \mu^{2}+\dfrac{\anmo}{\an} \mu + \dfrac{\anmt}{\an} \rb^{1/2} \; .
\edis
\item
If the radius of the disk centered at the origin in the theorems both in this and the previous section is larger than the Cauchy radius, then it 
will obviously contain all the zeros of the polynomial and the other disk(s) can be ignored. This can easily be detected by a simple substitution 
of the radius of that disk in equation~(\ref{sk_eq}).
\end{itemize}

\begin{figure}[H]
\begin{center}
\raisebox{0ex}{\includegraphics[width=0.475\linewidth]{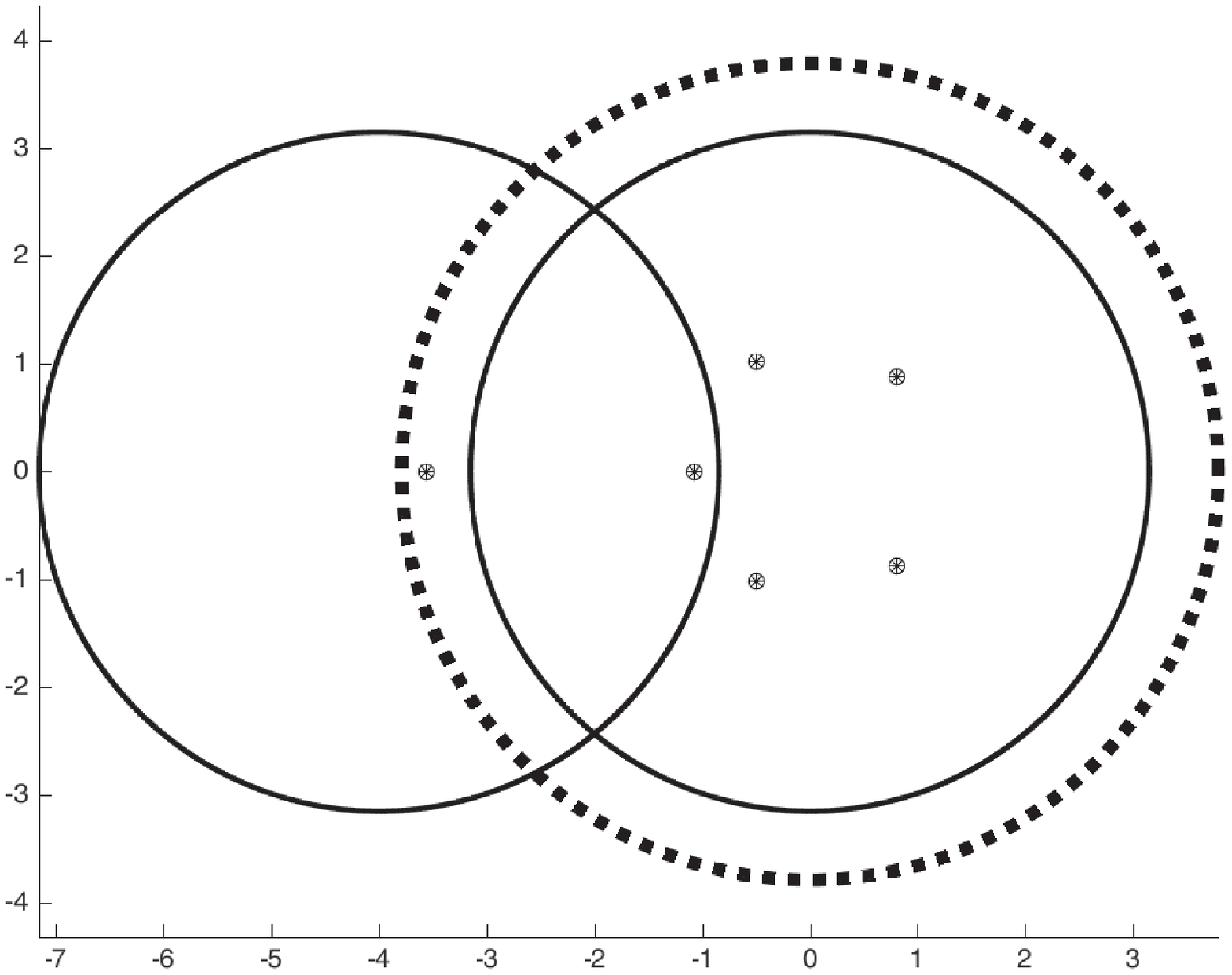}}
\hskip 0.5cm
\raisebox{0ex}{\includegraphics[width=0.475\linewidth]{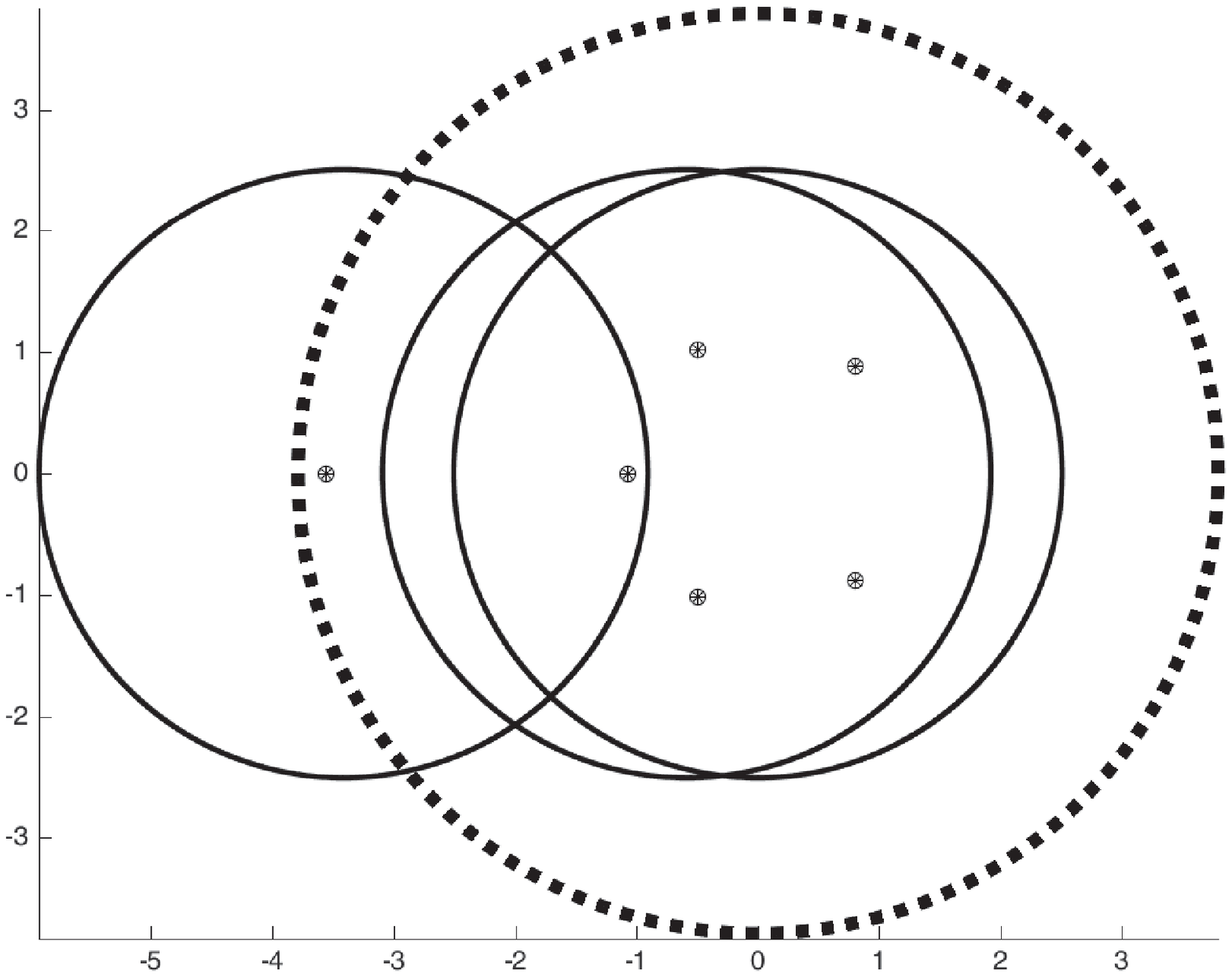}}
\caption{Comparison of Theorems~\ref{theorem_lemniscates1} and \ref{theorem_lemniscates4} for $p_{3}$.}
\label{eek_fig3}                      
\end{center}
\end{figure}

\begin{figure}[H]
\begin{center}
\raisebox{0ex}{\includegraphics[width=0.475\linewidth]{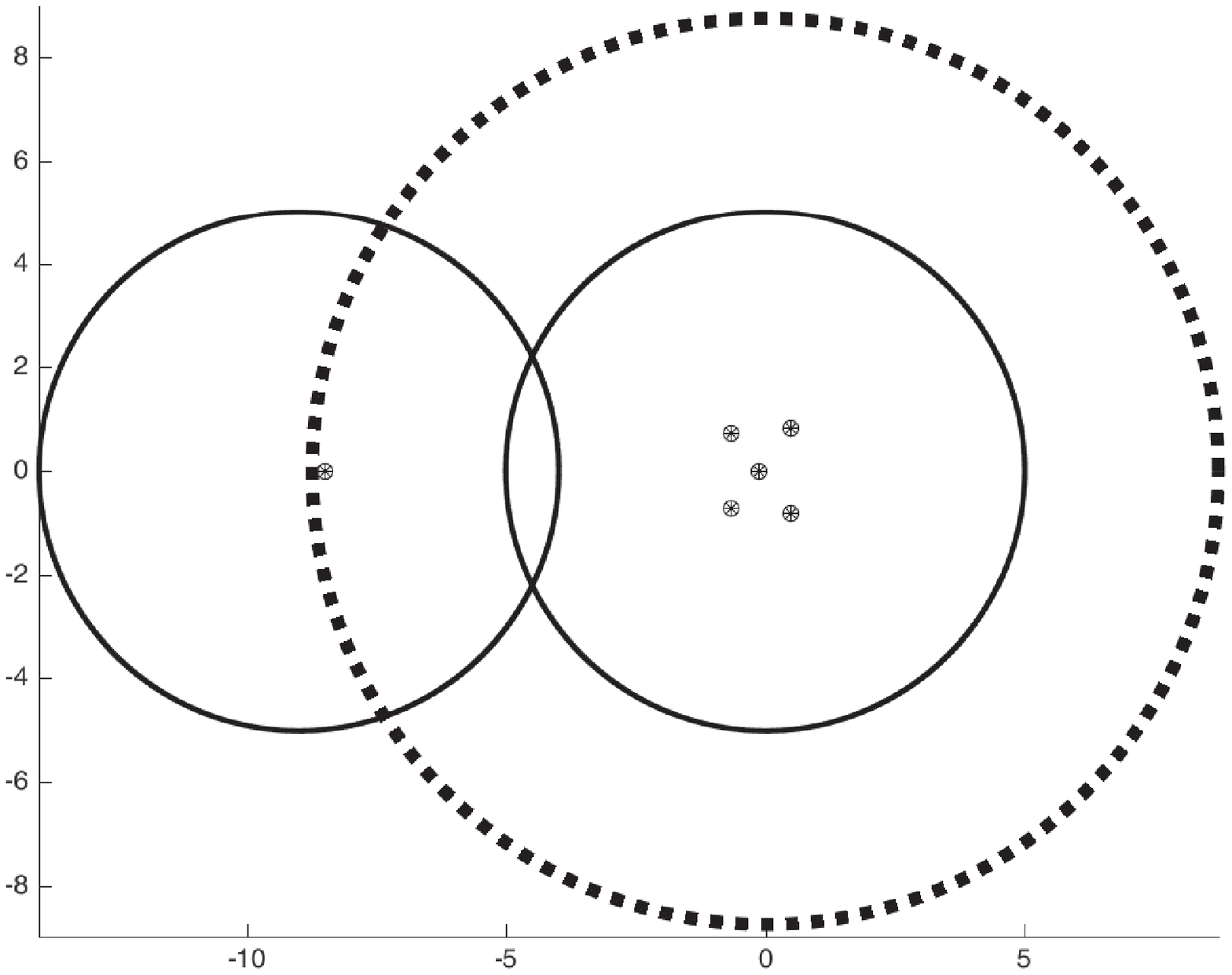}}
\hskip 0.5cm
\raisebox{0ex}{\includegraphics[width=0.475\linewidth]{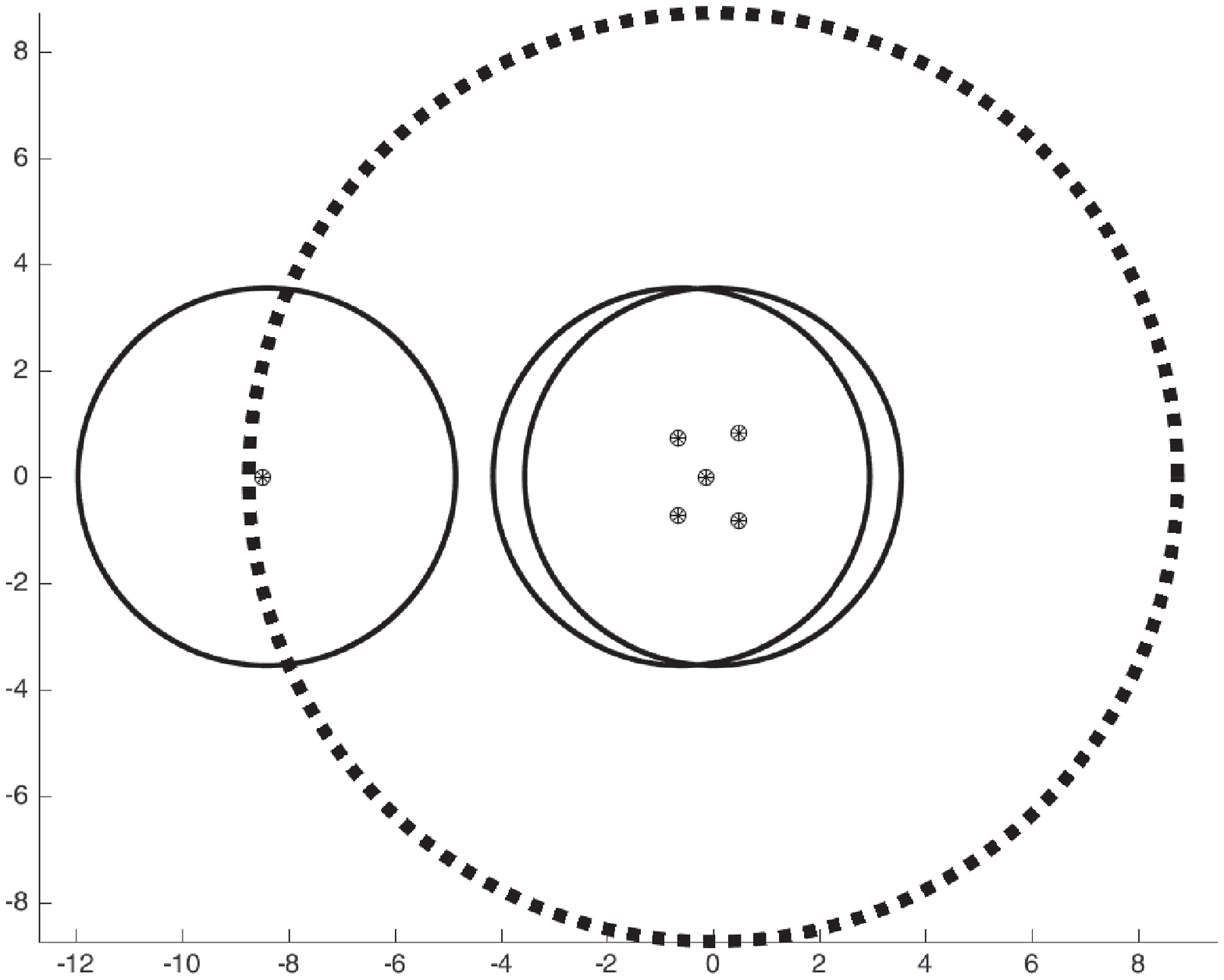}}
\caption{Comparison of Theorems~\ref{theorem_lemniscates1} and \ref{theorem_lemniscates4} for $p_{4}$.}
\label{eek_fig4}                      
\end{center}
\end{figure}

\vskip 0.25cm
\noindent {\bf Examples.} 
\begin{itemize}
\item
In~Figure~\ref{eek_fig3}, we compare the inclusion regions obtained from Theorem~\ref{theorem_lemniscates1} and Theorem~\ref{theorem_lemniscates4} 
for the same polynomial $p_{3}(z)=z^{6}+4z^{5}+2z^{4}+2z^{3}+3z^{2}+6z+7$ we used before. The dotted circle centered at the origin (with radius 3.788) 
is the boundary of the disk obtained from Theorem~\ref{Theorem_Bakchik2}(2). The solid circles in the figure on the left (with radii 3.151) are those 
obtained from Theorem~\ref{theorem_lemniscates1}, while those in the figure on the right (with radii 2.507) are obtained from 
Theorem~\ref{theorem_lemniscates4}. The zeros of $p$ are indicated by circled asterisks. The Cauchy radius in this case is 4.580.
\item
In~Figure~\ref{eek_fig4}, we carried out the same comparison as in Figure~\ref{eek_fig3} for the 
polynomial $p_{4}(z)=z^{6}+9z^{5}+5z^{4}+6z^{3}+4z^{2}+8z+1$. The dotted circle centered at the origin (with radius 8.744) 
is the boundary of the disk obtained from Theorem~\ref{Theorem_Bakchik2}(2). The solid circles in the figure on the left (with radii 5.012) are those 
obtained from Theorem~\ref{theorem_lemniscates1}, while those in the figure on the right (with radii 3.552) are obtained from
Theorem~\ref{theorem_lemniscates4}. 
As before, the zeros of $p$ are indicated by circled asterisks. The Cauchy radius of $p_{4}$ is 9.592. 
Here Theorem~\ref{theorem_lemniscates4} isolates one zero of $p_{4}$ from the others.
\end{itemize}

%
%

%
%

Theorem~\ref{theorem_lemniscates4} is just one more example of the kind of results that can be generated in the framework we have established,
which allows for many variations. More disks can be obtained by combining various multipliers with Theorem~\ref{theorem_bakchik_cauchy}.

\vskip 0.25cm \noindent {\bf Conclusion.} We have constructed a framework to derive generalizations of the classical \eneskakk theorem using two simple tools:
polynomial multipliers and a theorem establishing inclusion regions for the zeros of a polynomial. This framework unifies and simplifies the derivation of these
generalizations, obtaining new as well as old results in the process, while transparently showing how more such results can be generated. One feature of our
results, namely, zero inclusion regions consisting of more than one disk, is not found in any of the existing generalizations of this theorem.

\section{Appendix}
\label{appendix}

We state and prove a lemma that was used in Section~\ref{onedisk_notorigin}
\begin{lemma}
\label{lemma_recip}
The set $\lcb z \in \complex : \left | z + a  \right | \leq |a|R|z| \rcb$, 
where $a \in \complex$ and $R > 1/|a|$, is the closed exterior of a disk 
with center $a/(|a|^{2}R^{2}-1)$ and radius $|a|^{2}R/(|a|^{2}R^{2}-1)$. 
\end{lemma}
\prf
\begin{eqnarray*}
|z+a| \leq |a|R|z| 
& \Longleftrightarrow & |z+a|^{2} \leq |a|^{2}R^{2}|z|^{2}    \\ 
& \Longleftrightarrow & |z|^{2} + 2 Re(\bar{a}z) + |a|^{2} \leq |a|^{2}R^{2}|z|^{2}  \\ 
& \Longleftrightarrow & 0 \leq \lb |a|^{2}R^{2}-1 \rb |z|^{2} - 2 Re(\bar{a}z) - |a|^{2}     \\ 
& \Longleftrightarrow & 0 \leq |z|^{2} - \dfrac{2 Re(\bar{a}z)}{|a|^{2}R^{2}-1} - \dfrac{|a|^{2}}{|a|^{2}R^{2}-1}  \\ 
& \Longleftrightarrow & 0 \leq \left | z - \dfrac{a}{|a|^{2}R^{2}-1} \right |^{2} - \dfrac{|a|^{2}}{\lb |a|^{2}R^{2}-1 \rb^{2}}  
                                                                                                          - \dfrac{|a|^{2}}{|a|^{2}R^{2}-1}  \\ 
& \Longleftrightarrow & \left | z - \dfrac{a}{|a|^{2}R^{2}-1} \right |^{2} \geq \dfrac{|a|^{4}R^{2}}{\lb |a|^{2}R^{2}-1 \rb^{2}}  \\
& \Longleftrightarrow & \left | z - \dfrac{a}{|a|^{2}R^{2}-1} \right | \geq \dfrac{|a|^{2}R}{|a|^{2}R^{2}-1}  \; . \; \hskip 4cm \qed \\
\end{eqnarray*}

\end{document}